\documentclass{gtart}

\def\ifplaintex{\expandafter\ifx\csname documentclass\endcsname\relax}

\def\gtp{{\mathsurround=0pt\it $\cal G\mskip-2mu$eometry \&\ 
$\cal T\!\!$opology $\cal P\!$ublications}}  

\def\recd{{\small Received:\qua\receiveddate\ifx\reviseddate\relax
\else\qquad Revised:\qua\reviseddate\fi\par}} 

\def\lognumber#1{\def\thelognumber{#1}}
\def\volumenumber#1{\def\thevolumenumber{#1}}
\def\volumeyear#1{\def\thevolumeyear{#1}}
\def\papernumber#1{\def\thepapernumber{#1}}
\def\pagenumbers#1#2{\def\startpage{#1}\def\finishpage{#2}}
\def\published#1{\def\publishdate{#1}}

\def\received#1{\def\receiveddate{#1}}
\def\revised#1{\def\reviseddate{#1}}
\def\accepted#1{\def\accepteddate{#1}}

\long\def\asciiabstract#1{\long\def\theasciiabstract{#1}}
\def\asciikeywords#1{\def\theasciikeywords{#1}}

\let\\\par\let\thelognumber\relax\let\thevolumenumber\relax
\let\thepapernumber\relax\let\thevolumeyear\relax\let\startpage\relax
\let\finishpage\relax\let\publishdate\relax\let\receiveddate\relax
\let\reviseddate\relax\let\accepteddate\relax\let\theasciititle\relax
\let\theasciiauthors\relax
\let\theasciiabstract\relax\let\theasciikeywords\relax

\let\theasciiemail\relax

\ifplaintex
\font\logobig=cmssbx10 scaled 3836
\font\logomed=cmssbx10 scaled 2557
\else
\font\logobig=cmssbx10 scaled 4200
\font\logomed=cmssbx10 scaled 2800
\fi

\long\def\makeagttitle{   
\count0=\startpage
\agt\hfill      
\hbox to 45truept{\vbox to 0pt{\vglue -13truept{\logomed A\kern -.37em{\logobig 
T}\kern -.38em G}\vss}\hss}
\break
{\small Volume \thevolumenumber\ (\thevolumeyear)
\startpage--\finishpage\nl
Published: \publishdate}

\vglue .25truein

{\parskip=0pt\leftskip 0pt plus
1fil\def\\{\par\smallskip}{\Large\bf\thetitle}\par\medskip} \vglue
0.05truein

{\parskip=0pt\leftskip 0pt plus 1fil\def\\{\par}{\sc\theauthors}
\par\medskip}%
 
\vglue 0.03truein 

{\small\leftskip 25truept\rightskip 25truept{\bf Abstract}\stdspace\theabstract

{\bf AMS Classification}\stdspace\theprimaryclass
\ifx\thesecondaryclass\relax\else; \thesecondaryclass\fi\par
{\bf Keywords}\stdspace \thekeywords\par}\vglue 7truept

}   

\ifplaintex
\font\phead=cmsl9 scaled 950
\font\pnum=cmbx10 scaled 913
\font\pfoot=cmsl9 scaled 950
\headline{\vbox to 0pt{\vskip -4.5mm\line{\small\phead\ifnum
\count0=\startpage ISSN 1472-2739 (on-line) 1472-2747 (printed)
\hfill {\pnum\folio}\else\ifodd\count0\def\\{ }%
\ifx\theshorttitle\relax\thetitle\else\theshorttitle\fi\hfill{\pnum\folio}
\else\def\\{ and }{\pnum\folio}\hfill\ifx\theshortauthors\relax\theauthors
\else\theshortauthors\fi\fi\fi}\vss}}
\footline{\vbox to 0pt{\vglue 0mm\line{\small\pfoot\ifnum\count0=\startpage
\copyright\ \gtp\hfill\else
\agt, Volume \thevolumenumber\ (\thevolumeyear)\hfill\fi}\vss}}
\else
\font=cmsl9 scaled 1050
\font\lnum=cmbx10 
\font=cmsl9 scaled 1050
\makeatletter
\def\@oddhead{{\small\lhead\ifnum\count0=\startpage ISSN 1472-2739 
(on-line) 1472-2747 (printed)\hfill {\lnum\number\count0}\else\ifodd\count0
\def\\{ }\ifx\theshorttitle\relax \thetitle \else\theshorttitle\fi\hfill
{\lnum\number\count0}\else\def\\{ and }{\lnum\number\count0}
\hfill\ifx\theshortauthors\relax 
\theauthors\else\theshortauthors\fi\fi\fi}}\def\@evenhead{\@oddhead}
\def\@oddfoot{\small\lfoot\ifnum\count0=\startpage\copyright\ \gtp\hfill\else
\agt, Volume \thevolumenumber\ (\thevolumeyear)\hfill\fi}
\def\@evenfoot{\@oddfoot}
\makeatother
\fi
\let\maketitlepage\makeagttitle
\let\makeshorttitle\maketitlepage
\let\maketitle\maketitlepage

\newwrite\gtoutfile
\long\gdef\makeheadfile{  
{\def\\{, }\def\s{ }
\immediate\openout\gtoutfile head.xxx
\immediate\write\gtoutfile{To: math@arxiv.org}
\immediate\write\gtoutfile{Subject: put OR rep NNNNN:ppppp}
\immediate\write\gtoutfile{--text follows this line--}
\immediate\write\gtoutfile{Proxy-for: \ifx\theasciiauthors\relax
\theauthors\else\theasciiauthors\fi\s<\ifx\theasciiemail\relax\theemail\else\theasciiemail\fi>}
\immediate\write\gtoutfile{\noexpand\\}
\immediate\write\gtoutfile{Authors: \ifx\theasciiauthors\relax
\theauthors\else\theasciiauthors\fi}
{\def\\{ }\immediate\write\gtoutfile{Title: \ifx\theasciititle\relax
\thetitle\else\theasciititle\fi}}
\immediate\write\gtoutfile{Subj-class: GT or SG, GR etc}
\immediate\write\gtoutfile{MSC-class: \theprimaryclass\ifx\thesecondaryclass\relax\else, \thesecondaryclass\fi}
\immediate\write\gtoutfile{Journal-ref: Algebr. Geom. Topol. \thevolumenumber\s
(\thevolumeyear) \startpage-\finishpage}
\immediate\write\gtoutfile{Comments: Published by Algebraic and
Geometric Topology at}
\immediate\write\gtoutfile{\s\s\s  http://www.maths.warwick.ac.uk/agt/AGTVol\thevolumenumber/agt-\thevolumenumber-\thepapernumber.abs.html}
\immediate\write\gtoutfile{\noexpand\\}
\immediate\write\gtoutfile{}
\ifx\theasciiabstract\relax
\immediate\write\gtoutfile{\theabstract}\else
\immediate\write\gtoutfile{\theasciiabstract}\fi
\immediate\write\gtoutfile{}
\immediate\write\gtoutfile{\noexpand\\}
\immediate\write\gtoutfile{}
\immediate\closeout\gtoutfile}}  

\def\maketitlepage{\makeagttitle\makeheadfile}
\let\makeshorttitle\maketitlepage
\let\maketitle\maketitlepage

\lognumber{7}
\volumenumber{3}
\volumeyear{2003}
\papernumber{7}
\published{21 February 2003}
\pagenumbers{155}{185}
\received{24 November 2002}
\revised{27 January 2003}
\accepted{13 February 2003}

\usepackage{amssymb,amsmath,amscd}

\newcommand{\spin}{$Spin^c$-structure }
\newcommand{\spins}{$Spin^c$-structures }
\newcommand{\sqrr}{\sqrt{r}}
\newtheorem{theorem}{Theorem}[section]
\newtheorem{proposition}[theorem]{Proposition}
\newtheorem{lemma}[theorem]{Lemma}
\newtheorem{corollary}[theorem]{Corollary}
\newtheorem{assumption}[theorem]{Assumption}
\theoremstyle{remark}
\newtheorem{remark}[theorem]{Remark}

\begin{document}
\title{Grafting Seiberg-Witten monopoles}
\author{Stanislav Jabuka}
\address{Department of Mathematics, Columbia University\\2990 
Broadway, New York, NY 10027, USA}
\email{jabuka@math.columbia.edu}
\begin{abstract}
We demonstrate that the operation of taking disjoint unions of $J$-holomorphic curves (and thus
obtaining new 
$J$-holomorphic curves) has a Seiberg-Witten counterpart. The main theorem asserts that, given
two solutions 
$(A_i, \psi _i)$, $i=0,1$ of the Seiberg-Witten equations for the \spins $W^+ _{E_i}= E_i \oplus (E_i \otimes K^{-1})$
(with certain restrictions), 
there is a solution $(A, \psi )$ of the Seiberg-Witten equations for the \spin $W_E$ with $E= E_0\otimes E_1$, obtained by 
\lq\lq grafting\rq\rq  the two solutions  $(A_i, \psi _i)$. 
\end{abstract}

\asciiabstract{We demonstrate that the operation of taking disjoint
unions of J-holomorphic curves (and thus obtaining new J-holomorphic
curves) has a Seiberg-Witten counterpart. The main theorem asserts
that, given two solutions (A_i, psi _i), i=0,1 of the Seiberg-Witten
equations for the Spin^c-structure W^+_{E_i}= E_i direct sum (E_i
tensor K^{-1}) (with certain restrictions), there is a solution (A,
psi) of the Seiberg-Witten equations for the Spin^c-structure W_E with
E= E_0 tensor E_1, obtained by `grafting' the two solutions (A_i,
psi_i).}

\asciikeywords{Symplectic 4-manifolds, Seiberg-Witten gauge theory,
J-holomorphic curves}

\maketitle
\let\\\par
\section{Introduction} \label{intro}

In his series of groundbreaking works \cite{kn:taub1},
\cite{kn:taub2}, \cite{kn:taub3}, Taubes showed that the
Seiberg-Witten invariants and the Gromov-Witten invariants (as defined
in \cite{kn:taub4}) for a symplectic 4-manifold $(X, \omega)$ are the
same. His results opened the door to a whole new world of interactions
between the two theories that had previously only been
speculations. The most spectacular outcomes of this interplay were new
results that in one theory were obvious but when translated into the
other theory, became highly nontrivial. An example of such a
phenomenon is the simple formula relating the Seiberg-Witten invariant
of a \spin $W$ to the Seiberg-Witten invariant of its dual \spin
$W^*$, i.e.\ the one with $c_1 (W^*) = -c_1 (W)$. The formula reads:
$$ SW_X (W^*) = \pm \,SW_X(W) $$
When translated into the Gromov-Witten language, this duality becomes 
\begin{equation}
Gr_X (E) = \pm \, Gr_X (K-E) \label{eq:grdual}
\end{equation}
Here $K$ is the canonical class of $(X,\omega )$ and $E\in H^2 (X;\mathbb{Z})$ is related to
$W$ as $c_1(W^+)=2\,E - K$. This is a highly nonobvious result about $J$-holomorphic curves, 
even in the simplest case when $E=0$. In that case we obtain that $Gr_X (K) = \pm \, Gr _X (0)
=\pm \, 1$, the latter equation simply being the definition of 
$Gr_X (0)$. This gives an existence result of a $J$-holomorphic representative for the class $K$, a
result unknown prior to Taubes' theorem. The formula 
\eqref{eq:grdual} has recently been proved by S. Donaldson and I. Smith \cite{kn:donaldson}
without any reference to Seiberg-Witten theory (but under slightly stronger
restrictions on $(X, \omega )$ than in Taubes' theorem). 

In the author's opinion, proving a result about Gromov-Witten theory which had only been
known through its relation with Seiberg-Witten theory, without relying 
on the latter, has a number of benefits.  
One is to understand Gromov-Witten theory from within better. But also to possibly generalize
the theorem to a broader class of manifolds. Recall that Taubes'
theorem equates the two invariants only on symplectic 4-manifolds. Both
Seiberg-Witten and Gromov-Witten theory are defined over larger sets 
of manifolds, namely all smooth 4-manifolds and all symplectic manifolds (of any dimension)
respectively. On the other hand, even within the category of symplectic 
4-manifolds, one can hope for more nonvanishing theorems i.e.\ theorems of the type
$Gr_X(E)\ne 0$ for classes $E\ne 0,\, K$. 
The techniques used by Donaldson and Smith are promising in that direction.

The aim of this paper is to prove a result in the same vein but going the opposite direction.
Namely, on the Gromov-Witten side, given two classes 
$E_i \in H^2(X;\mathbb{Z} )$, $i=0, 1$ with  $E_0 \cdot E_1= 0 $ and $J$-holomorphic curves
$\Sigma _i$ with $[\Sigma _i ] = P. D. (E_i)$,
one can define a new $J$-holomorphic curve $\Sigma = \Sigma_0 \sqcup \Sigma _1$. By the
assumption $E_0\cdot E_1 = 0$, the two curves $\Sigma _i$ 
are either disjoint or share toroidal components (see \cite{kn:mcduff}). In the former case,
$\Sigma$ is simply the disjoint union of $\Sigma _0$ and $\Sigma _1$ and in the 
latter case one needs to  replace the  tori shared by $\Sigma _0$ and $\Sigma _1$ with their
appropriate multiple covers.
This induces a map on moduli spaces: 
\begin{equation}
{\cal M}^{Gr} _X (E_0)\times {\cal M}^{Gr} _X (E_1) \stackrel{\sqcup}{\rightarrow} {\cal
M}^{Gr} _X (E_0 + E_1) \label{eq:grmap}
\end{equation} 
This article describes the Seiberg-Witten counterpart of \eqref{eq:grmap}. That is, given two
complex line bundles 
$E_0$ and $E_1$ (with certain restrictions, see the assumption \ref{assump:mainassump}  below for a precise
statement) and two solutions $(A_i, \psi _i)$ of the Seiberg-Witten equations 
for the \spins $W_{E_i} = E_i \oplus (E_i \otimes K^{-1} )$, $i=0,1$ and with Taubes' large $r$
perturbation, we show how to produce a solution $(A,\psi) = 
(A_0, \psi _0)\cdot (A_1, \psi _1)$ 
for the \spin $W_E$ with $E=E_0\otimes E_1$. This operation induces the following
commutative diagram:\eject
\noindent\hbox{}\vspace{-10pt}
\begin{equation}
\begin{CD} 
 {\cal M}_X ^{SW} (E_0) \times {\cal M}_X^{SW} (E_1) @>{\cdot}>> {\cal M}_X^{SW}
(E_0\otimes E_1)  \\
@V{\Theta}VV                                                
          @VV{\Theta}V \\
 {\cal M}_X ^{Gr} (E_0) \times {\cal M}_X^{Gr} (E_1) @>{\sqcup}>> {\cal M}_X^{Gr}
(E_0 +  E_1)  \\
\end{CD}  \label{eq:grsw}
\end{equation}

Here the map $\Theta:{\cal M}_X^{SW}(E) \rightarrow {\cal M}_X^{Gr}(E)$ is the map
described in \cite{kn:taub1} that associates to each solution of the 
Seiberg-Witten equations an embedded $J$-holomorphic curve. The solution $(A,\psi )$ is
obtained by ``grafting" the two solutions $(A_i,\psi _i )$. The key observation here is that for the
large $r$ version of Taubes'
perturbation, a solution $(B,\phi )$ of the Seiberg-Witten equations for the \spin $W_E$ is
``concentrated" near the zero set of 
$\sqrr \, \alpha$, the $E$ component of $\phi$. That is, the restriction of $(B,\phi )$ to the
complement of a regular neighborhood of $\alpha ^{-1} (0)$ 
converges pointwise (under certain bundle identifications) to the unique solution $({\cal A}_0,
\sqrr \, u_0)$ for the anticanonical 
\spin $W_0=\underline{\mathbb{C}} \oplus 
K^{-1}$. This is used to define a first approximation of $\psi$ by declaring it to be equal to $\psi
_i$ in a regular neighborhood $V_i$ of $\alpha _i ^{-1} (0)$ and 
equal to $\sqrr \, u_0$ on the complement of $V_0 \cup V_1$. Bump functions are used to
produce a smooth spinor. The first approximation of $A$ is simply the 
product connection $A_0 \otimes A_1$. The contraction mapping principle is then evoked to
deform this approximate solution to an
 honest solution of the 
Seiberg-Witten equations. 
The author has learned the techniques employed in this article  from the inspiring
work of Taubes on gauge theory of symplectic 4-manifolds, most notably from \cite{kn:taub2}. 

The article is organized as follows. In section 2 we review the needed Seiberg-Witten theory on
symplectic 4-manifolds. Section 3 explains how to define 
an ``almost" monopole $(A',\psi ')$ from a pair of monopoles $(A_i,\psi_i)$, $i=0,1$. It also
analyzes the asymptotic (as $r\rightarrow \infty$) regularity theory 
for the linearized operators $L_{(A_i,\psi _i)}$ and deduces a corresponding result for
$L_{(A',\psi ' )}$. The latter is used in combination with the contraction
mapping principle to obtain an ``honest" monopole $(A, \psi )$.  Section 4 compares the present method of
grafting monopoles to the one used in exploring Seiberg-Witten 
theory on manifolds $X$ which are obtained as a fiber sum: $X=X_1 \# _{\Sigma} X_2$.
Section 5 proves a converse to theorem \ref{theorem:main}. It explains which monopoles 
in the \spin $W_E$ can be obtained as products of monopoles $(A_i,\psi _i)$ in the \spins
$W_{E_i}$, $i=0,1$ with $E_0\otimes E_1 = E$ and with the property that $\Theta (A_i, \psi _i)$ does
not contain multiply covered tori. 

\medskip
{\bf Acknowledgment}\qua The author would like to express his gratitude to his thesis advisor,
professor Ron Fintushel, for his continuing help and 
encouragement during the process of writing this article, the author's doctoral thesis. 

\section{The Seiberg-Witten equations on symplectic manifolds}  \label{sec:1}
\subsection{Introduction} \label{sec:21}
Let $(X,\omega )$ be a symplectic, smooth, compact 4-manifold with symplectic form $\omega
$. Denote by ${\cal J}$ the set of all 
almost complex structures $J$ on $TX$ that are compatible with $\omega$, i.e.\ the ones for
which 
$$g(v,w) = g_J (v,w)= \omega (v, Jw )\quad \quad v,w\in TX$$ 
defines a Riemannian metric on $X$.  Given a $J\in {\cal J}$, the associated metric $g_J$ will
always be assumed throughout to be the metric of choice. 

On any almost complex 4-manifold there is a anticanonical \spin $W_0=W_0 ^+ \oplus W_0 ^-$
determined by the almost complex structure as:
\begin{align}
W_0 ^+ = & \, \Lambda ^{0,0} \oplus \Lambda ^{0,2} = \underline{\mathbb{C}} \oplus K^{-1} 
 \nonumber \\
W_0 ^- = & \,  \Lambda ^{0,1} \nonumber \\
v.\alpha =  &  \sqrt{2} \,\left(  v^* _{0,1} \wedge \, \alpha - \iota _v \alpha \right)
\quad \quad v\in T_x X, \, \alpha \in W_{0,x}, \, x\in X \nonumber 
\end{align}
In the above, $v^* _{0,1} \in \Lambda ^{0,1}$ denotes the (0,1) projection of $v^*\in T_x ^*
X$, the dual of $v\in T_xX$.
All other \spins can be obtained from $W_0$ by tensoring it with a complex line bundle $E$ and
extending Clifford multiplication trivially over the $E$
factor, i.e. 
\begin{align}
W_E ^{\pm} = & \, E\otimes W_0 ^{\pm} \nonumber \\
v.(\varphi \otimes \alpha) = & \,  \varphi \otimes (v.\alpha) \quad \quad \varphi \in E_x, \, v\in
T_x X,\, \alpha \in W_{0,x}, \, x\in X \nonumber 
\end{align}
The symplectic form $\omega$ induces a splitting of $\Lambda ^{2,+}$ as 
\begin{equation}
 \Lambda ^{2,+} \cong \mathbb{R}\cdot \omega \oplus \Lambda ^{0,2} \label{eq:curvspli}
\end{equation}
which will be used below to write the curvature component of the Seiberg-Witten equations as
two equations, one for each of the summands on the right-hand side 
of \eqref{eq:curvspli}. 

Given a \spin $W_E$ on $X$, the Seiberg-Witten equations are a coupled, elliptic  system of
equations for a pair $(A,\psi )$ of a connection $A$ on $E$ and
a positive spinor $\psi \in \Gamma (W^+ _E) = \Gamma (E\oplus (E\otimes K^{-1}))$. 
The connection $A$ on $E$ together with a fixed connection ${\cal A}_0$ on $K^{-1}$ (which
will be made specific in a bit), induces a $Spin^c$-connection on $W_E$ which we will denote 
by $\nabla ^A$ and which in turn gives rise to the Dirac operator $D_A : \Gamma (W_E)
\rightarrow \Gamma (W_E)$. It proves convenient to write the spinor $\psi$ in the form 
$$\psi = \sqrr (\alpha , \beta ) \quad \quad \alpha \in \Gamma (E), \, \beta \in \Gamma (E\otimes
K^{-1}) $$
where $r\ge 1$ is a parameter whose significance will become clear later. With this understood,
the Seiberg-Witten equations read:
\begin{align}
D_A (\psi ) =  & \,0 \nonumber  \\
F_A ^+ =  & \, q(\psi , \psi ) + \mu \label{eq:sw}
\end{align}
Here $\mu \in i\Omega ^{2,+}$ is a fixed imaginary, self-dual two form on $X$ and $q:\Gamma
(W_E^+) \times \Gamma (W_E^+) \rightarrow i\Omega ^{2,+}$
is the bilinear quadratic map given explicitly  by
\begin{equation}
q(\psi , \psi ) = \frac{ir}{8} (|\alpha |^2 - |\beta |^2 )\omega + \frac{ir}{4} (\bar{\alpha } \beta + 
\alpha \bar{\beta} )
\end{equation}
\subsection{The anticanonical $Spin^c$-structure} \label{sec:22}
It is another result of Taubes' \cite{kn:taub5} that the Seiberg-Witten invariant of the
anticanonical \spin on a symplectic manifold is equal to $\pm 1$.
Furthermore, the equations have exactly one solution $({\cal A}_0,\sqrr\cdot u_0)$, $u_0 \in
\Gamma (\underline{\mathbb{C}} )$,  for the choice of 
\begin{equation}
\mu=F_{{\cal A}_0}^+-\frac{i\,r}{8} \, \omega \label{eq:per}
\end{equation}
in \eqref{eq:sw} and for $r\gg 1$. The 
purpose of this section is to describe the solution $({\cal A}_0, \sqrr \, u_0)$ and its linearized
operator. 

The pair $({\cal A}_0 , \sqrr\cdot u_0)$ is characterized (up to gauge) by the condition 
\begin{equation}
\langle \nabla ^0 u_0 , u_0\rangle = 0  \label{eq:canon}
\end{equation}
 (where $\nabla ^0$ is the 
$Spin^c$-connection induced by ${\cal A}_0$) and can be obtained as follows: 
let $u_0$ be any section of $\underline{\mathbb{C}}\oplus K^{-1}$ with $|u_0| = 1$ and 
whose projection onto the second summand is zero. Likewise, let $A$ be any connection on
$K^{-1}$ and let $\nabla ^A$ be its induced 
$Spin^c$-connection on $W^+ _0 =  \underline{\mathbb{C}}\oplus K^{-1}$. Set $a=\langle
u_0 , \nabla ^A u_0 \rangle$. 
This defines an imaginary valued 1-form as can easily be seen: 
$$ a + \bar{a} = \langle \nabla ^A u_0 ,  u_0 \rangle + \langle u_0 , \nabla ^A u_0 \rangle  = d\,
|u_0|^2 = 0 $$
Define the connection ${\cal A}_0$ on $K^{-1}$ by ${\cal A}_0 = A - a$ which induces the
$Spin^c$-connection $\nabla ^0 = \nabla ^A -a $ on $W^+ _0$. 
This connection clearly satisfies \eqref{eq:canon}. With the choice of $\mu $ as in
\eqref{eq:per}, the Seiberg-Witten equations \eqref{eq:sw} take the form: 
\begin{align}
D_A \psi = & \, 0 \nonumber \\
F_A ^+ = & \, \frac{ir}{8} (|\alpha| ^2 -1 - |\beta |^2) \omega + F_{{\cal A}_0}^+ + \frac{ir}{4}
(\bar{ \alpha} \beta + \alpha \bar{\beta } ) \label{eq:swsymp}
\end{align} 
Since the $\beta$-component of $u_0$ is zero and since $|\alpha| = |u_0| = 1$, the pair $({\cal
A}_0 , u_0 )$ clearly solves the second equation of \eqref{eq:swsymp}.
The fact that is also solves the first equation relies on the property $d \, \omega = 0$ of  $\omega$ as well as
\eqref{eq:canon}. 
Taubes \cite{kn:taub5} showed that there are, up to gauge, no other solutions to
\eqref{eq:swsymp} and, as we shall presently see, that the solution 
$({\cal A}_0 , u_0)$ is a smooth 
solution in the sense that the linearization of \eqref{eq:swsymp} at $({\cal A}_0 , u_0 )$ has
trivial cokernel. These two facts together show that $SW_X (W_0) = \pm 1$.

Define $S:L^{1,2} (i\Lambda ^1 \oplus W^+ _0) \rightarrow L^2 (i\Lambda ^0 \oplus i\Lambda
^{2,+} \oplus W^- _0 )$ to be the linearized Seiberg-Witten 
operator for the solution $({\cal A}_0, u_0)$. Thus, for $(b,(\xi _0 , \xi _2))\in L^{1,2}
(i\Lambda ^1 \oplus (\underline{\mathbb{C}}\oplus K^{-1})) $ we have:
\begin{align}  \label{linear:easy}   
                                & \quad \left(  d^* b + i \frac{\sqrr}{\sqrt{2}} \mbox{ Im} (\bar{u}_0 \xi _0 ) , \right.   \nonumber  \\
S(b,(\xi _0 , \xi _2)) = & \quad d^+ b - \sqrr q(\xi , u_0) - \sqrr q( u_0 , \xi ) ,  \\ 
                                & \quad  \left. D_{A_0} (\xi _0 , \xi _2) + \frac{\sqrr}{2} \, b.u_0  \right)   \nonumber  
\end{align}
Let $S^* :  L^2 (i\Lambda ^0 \oplus i\Lambda ^{2,+} \oplus W^- _0) \rightarrow L^{1,2}
(i\Lambda ^1 \oplus W^+ _0)$ be the formal adjoint of $S$. 
The following proposition and corollary are proved in \cite{kn:taub2}, section 4. 
\begin{proposition} Let $S$ and $S^*$ be as above. Then the operator $SS^*$ on $L^2 (i\Lambda ^0
\oplus i\Lambda ^{2,+} \oplus W^- _0)$ is given by 
\begin{equation}
 SS^* = \frac{1}{4} \nabla ^{0,*} \nabla ^0 + {\cal R}_0 + \sqrr {\cal R}_1 + \frac{r}{8}
\label{eq:laplac}
\end{equation}
where $\nabla ^{0,*}$ is the adjoint of $\nabla ^0$ and where ${\cal R}_i,\, i=0,1$ are
certain $r$-indep\-end\-ent endomorphism on 
$L^2 (i(\Lambda ^0 \oplus \Lambda ^{2,+}) \oplus W^- _0)$. 
\end{proposition}

The proof  is a straightforward calculation, terms of the form $D_{{\cal A}_0} D_{{\cal A}_0}
^*$ are simplified using the
Weitzenb\"ock formula for the Dirac operator. An important consequence of \eqref{eq:laplac} is
the following:
\begin{corollary}
With $S$ and $S^*$ as above, the smallest eigenvalue $\lambda _1$ of $SS^*$ is bounded from
below by $r/16$. In particular, $S$ is invertible
and $S^{-1}$ satisfies the bounds 
\begin{equation}
|| S^{-1} y ||_2  \le \frac{4}{\sqrr} \, ||y||_2 \quad \mbox{ and } \quad || S^{-1} y ||_{1,2}  \le C \,
||y||_2 \label{eq:estim}
\end{equation}  
where $C$ is $r$-independent.
\end{corollary}
\subsection{The general case and bounds on $(a,\psi )$} \label{sec:23}
Consider now a \spin $W_E = E\otimes W_0$ on $X$.  
The connection ${\cal A}_0$ on $K^{-1}$ and a choice of a connection  $B_0$ on $E$ together
induce a connection 
$B_0^{\otimes 2} \otimes {\cal A}_0$ on $E^{\otimes 2}\otimes K^{-1} = c_1 (W_E^+ )$ by
the product rule:
$$B_0^{\otimes 2} \otimes {\cal A}_0 (\varphi _1 \otimes \varphi _2 \otimes \phi ) = B_0
(\varphi _1) \otimes \varphi _2 \otimes \phi + 
                                                       \varphi _1 \otimes
B_0(\varphi _2 ) \otimes \phi + 
                                                  \varphi _1 \otimes
\varphi _2 \otimes {\cal A}_0(\phi ) $$
The space of connections on $E^{\otimes 2}\otimes K^{-1}$ is an affine space with associated
vector space $i\Omega ^1_X$. With the choice of a ``base" 
connection $B_0^{\otimes 2} \otimes {\cal A}_0$ in place, we will from now on regard
solutions to the Seiberg-Witten equations as pairs $(a,\psi )\in i \Omega ^1 _X 
\times \Gamma (W^+_ E)$ rather than $(A ,\psi ) \in \mbox{Conn}(E^{\otimes 2}\otimes
K^{-1}) \times \Gamma (W^+_ E)$, the relation between the two being: 
$$ A = B_0^{\otimes 2} \otimes {\cal A}_0 + a $$
We will agree to use henceforth the choice of $\mu $ in  \eqref{eq:sw} to be: 
\begin{equation}
\mu = - \frac{i r }{8} \omega + F_{{\cal A}_0} ^+  \label{eq:pergen} 
\end{equation}
For $\psi \in \Gamma (E \otimes (\underline{\mathbb{C}} \oplus K^{-1} ))$  we will write $\psi
= \sqrr \, (\alpha \otimes u_0 , \beta )$ with 
$\alpha \in \Gamma (E)$ and $\beta \in \Gamma (E\otimes K^{-1} )$ and $u_0$ as in the
previous section. 

With these conventions understood and with the use of \eqref{eq:curvspli}, the Seiberg-Witten
equations \eqref{eq:sw} become:
\begin{align}
D_a \psi = & \, 0 \nonumber \\
F_a^{1,1} = & \,  \frac{ir}{8} (|\alpha |^2 - |\beta |^2 -1 ) \, \omega  \label{eq:swgen}  \\ 
F_a^{0,2} = & \, \frac{ir}{4} \bar{\alpha } \beta \nonumber 
\end{align}  
Here $F_a^{i,j}$ is the orthogonal projection of $2\,F^+_{B_0} + d^+ a$ onto $\Lambda
^{i,j}$. The linearized operator $L_{(a,\psi)}:L^{1,2} (i\Lambda ^1 \oplus W^+ _E) \rightarrow L^2 (i\Lambda ^0 \oplus i\Lambda
^{2,+} \oplus W^- _E )$ of the Seiberg-Witten equations for a solution $(a,\psi)$ of \eqref{eq:swgen} is: 
\begin{align}  \label{linear:hard}   
                                                  & \quad \left(  d^* b + i \frac{\sqrr}{\sqrt{2}} \mbox{ Im} (\langle \alpha , \xi _0 \rangle + \langle \beta , \xi_2  \rangle ) , \right.   \nonumber  \\
L_{(a,\psi)}\, (b,(\xi _0 , \xi _2)) = & \quad d^+ b - \sqrr q(\xi_0 + \xi _2  , \psi) - \sqrr q( \psi , \xi_0 + \xi _2 ) ,  \\ 
                                                  & \quad  \left. D_{a} (\xi _0 + \xi _2) + \frac{\sqrr}{2} \, b.\psi  \right)   \nonumber  
\end{align}
It is another result of Taubes' that the operator $L_{(a,\psi)}$ has Fredholm index zero on a symplectic manifold with $b^+_2\ge 2$, provided that $E$ is a basic class. 
As we will use this fact repeatedly throughout the paper, we give a short proof of it here:

\begin{theorem}[Taubes] Let $X$ be a symplectic manifold with $b^+_2\ge 2$ and $E\in H^2(X;\mathbb{Z})$ a basic class, i.e.\ $SW(W_E)\ne 0$. Let $(a,\psi)$ be a
solution of \eqref{eq:swgen} and $L_{(a,\psi)}$ be the operator defined by \eqref{linear:hard}. Then the Fredholm index of $L_{(a,\psi)}$ is equal to zero. 
\end{theorem} 

\proof
As $E$ is assumed to be a Seiberg-Witten basic class, it has to also be a Gromov-Witten basic class. In particular, the dimension of the Gromov-Witten moduli space
has to be non-negative:
\begin{equation} \nonumber
\mbox{dim }{\cal M}^{Gr}(E) = \frac{1}{2} (E^2-K\cdot E) \ge 0  
\end{equation}
Let $\Sigma$ be an embedded $J$-holomorphic curve in $X$ with $[\Sigma] = P.D.(E)$. Then the adjunction formula for $\Sigma$ states:
\begin{equation}  \nonumber
2g-2=E^2+K\cdot E
\end{equation}
Combining these last two relations we obtain two inequalities:
\begin{equation} \nonumber
E^2\ge g-1 \quad \quad \mbox{and} \quad \quad K\cdot E \le g-1 
\end{equation}
Let $n\ge 0$ be the integer such that $E^2 = g-1+n$ and $K\cdot E = g-1-n$. Since $E$ is a Gromov-Witten basic class, by duality, so is $K-E$. But then 
(by positivity of intersection of $J$-holomorphic curves) we must have: 
$$0\le E\cdot (K-E) = E\cdot K - E^2 = g-1-n - (g-1+n) = -2n \le 0 $$ 
This forces $n=0$ and so $E^2=g-1 = K\cdot E$. Using these in the index formula for $L_{(A,\psi)}$ immediately yields the desired result:
$$ \mbox{Ind } L_{(a,\psi)}  = \frac{1}{4} \left( (2E-K)^2 - (3\sigma + 2e)\right) =  \frac{1}{4} \left( K^2 - (3\sigma + 2e)\right)  = 0 \eqno{\qed}$$

We also use this section to remind the reader of several useful bounds that a solution $(a,\psi )$
of the Seiberg-Witten equations satisfies. 
These bounds are provided courtesy of \cite{kn:taub1} and their proofs rely solely on properties
of the Seiberg-Witten equations.

A solution $(a, \psi )$ of \eqref{eq:swgen} satisfies the following bounds: 
\begin{align}
|\alpha | \le & \, 1 + \frac{C}{r}   \nonumber  \\
|\beta| ^2 \le & \, \frac{C}{r} (1-|\alpha|^2) + \frac{C'}{r^3}  \label{eq:mainbounds} \\
|\nabla ^A \alpha |^2_x \le & \, C\, \sqrr \, \exp \left(  -\frac{\sqrr}{C} \mbox{dist}(x,\alpha ^{-1}
(0)) \right) , \quad x\in X  \nonumber  \\
|1-|\alpha (x) |^2 | \le & \, C \exp \left( -\frac{\sqrr}{C} \mbox{dist}(x,\alpha ^{-1} (0)) \right),
\quad x\in X  \nonumber \\
|F_a| \le & \frac{r}{4\sqrt{2}} (1-|\alpha|^2) + C \nonumber \\
|F_a(x)| \le & C\, r \exp \left( -\frac{1}{C} \sqrt{r} \, \mbox{dist}(x, \alpha ^{-1} (0) ) \right), \quad x\in X  \nonumber
\end{align}
\begin{remark} The constant $C$ appearing above may change its value from line to line. It is  important  to point out that 
$C$  only depends on the \spin $W_E$ and the Riemannian metric
$g$ but {\bf not} on the particular choice of the parameter $r$. This will be the case for all the numerous constants (all labeled $C$) appearing 
subsequently and we will henceforth tacitly adopt this misuse of notation.  
\end{remark}
\section {The main part} \label{sec:3}
\subsection {Producing the approximate solution $(a,\psi )$ from a pair\nl $(a_0,\psi _0 )$, $(a_1,\psi _1)$ } \label{sec:31}
Let $E_0$ and $E_1$  be two complex line bundles over $X$. 
The aim of this section is to produce an approximate solution $(a,\psi)$ of the Seiberg-Witten
equations for the \spin $W_{E_0 \otimes E_1}$ 
from two solutions $(a_0,\psi _0)$ and $(a_1,\psi _1)$ for the \spins  $W_{E_0}$ and $W_{E_1}$
respectively. Implicit to our discussion are the choices of 
two ``base" connections $B_0$ and $B_1$ on $E_0$ and $E_1$ and the product connection
$B_0\otimes B_1$ they determine 
on $E_0\otimes E_1$. As before, we will write $\psi _i = 
\sqrt{r} (\alpha _i \otimes u_0 , \beta _i)$, $i=0$, $1$, and $\psi = \sqrt{r} (\alpha \otimes u_0,
\beta)$. We define $(a ,\psi)$ as:
\begin{align}
               a &= a_0 + a_1      \nonumber     \\
       \alpha & = \alpha _0 \otimes \alpha _1         \label{eq:def1} \\
         \beta & = \alpha _0 \otimes \beta _1 + \alpha _1 \otimes \beta _0   \nonumber  
\end{align}

The first task at hand is to  check how close $(a,\psi )$ comes to solving the Seiberg-Witten
equations. We begin by calculating $D_a \psi$  locally 
at a point $x\in X$. Choose an orthonormal frame $\{ e_i \}_i $ in a neighborhood of $x$ and let $\{
e^i \} _i$ be its dual frame. 
\begin{align}
     D_a (\psi)  & = \sqrr D_a (\alpha _0 \otimes \alpha _1 \otimes u_0 + \alpha _0 \otimes \beta _1
+ \alpha _1 \otimes \beta _0 )  \nonumber \\ 
                     & = \sqrr D _{a_0} (\alpha _0 \otimes u_0) \otimes \alpha _1 + \sqrr \alpha _ 0
\otimes D_{a_1} (\alpha _1 \otimes u_0) + \nonumber \\
                     & \phantom{ = } + \sqrr e^i .\nabla ^a _{e_i} (\alpha _0 \otimes \beta _1 + \alpha _1
\otimes \beta _0)  \nonumber \\ 
                     & = \sqrr D _{a_0} (\alpha _0 \otimes u_0) \otimes \alpha _1 + \sqrr \alpha _ 0
\otimes D_{a_1} (\alpha _1 \otimes u_0) + \nonumber \\
                     & \phantom{ = } + \sqrr (\alpha _0 \otimes e^i.(\nabla ^{a_1} _{e_i} \beta _1) +
\alpha _1 \otimes e^i.(\nabla ^{a_0} _{e_i} \beta _0) + \nonumber \\
                     & \phantom{ = } + (\nabla ^{a_0} _{e_i} \alpha _0)\otimes e^i .\beta_1 +  (\nabla
^{a_1} _{e_i} \alpha _1)\otimes e^i .\beta_0) \nonumber  \\
                     & = \sqrr D _{a_0} (\alpha _0 \otimes u_0) \otimes \alpha _1 + \sqrr \alpha _ 0
\otimes D_{a_1} (\alpha _1 \otimes u_0) + \nonumber \\
                     & \phantom{ = } + \sqrr (\alpha _0 \otimes D_{a_1} \beta _1 + \alpha _1 \otimes
D_{a_0} \beta _0 ) + \nonumber \\ 
                     & \phantom{ = } + \sqrr ( (\nabla ^{a_0} _{e_i} \alpha _0)\otimes e^i .\beta_1 + 
(\nabla ^{a_1} _{e_i} \alpha _1)\otimes e^i .\beta_0) \nonumber  \\ 
                     & = (D_{a_0} \psi _0 )\otimes \alpha _1 + \alpha _0 \otimes (D_{a_1} \psi _1) +
\nonumber \\ 
                     & \phantom{ = } + \sqrr ( (\nabla ^{a_0} _{e_i} \alpha _0)\otimes e^i .\beta_1 + 
(\nabla ^{a_1} _{e_i} \alpha _1)\otimes e^i .\beta_0) \nonumber   \\
                     & = \sqrr  (\nabla ^{a_0} _{e_i} \alpha _0)\otimes e^i .\beta_1 + \sqrr (\nabla ^{a_1}
_{e_i} \alpha _1)\otimes e^i .\beta_0  \label{eq:dirac-point} 
\end{align}

It is easy to see, using the bounds in \eqref{eq:mainbounds}, that the first term in
\eqref{eq:dirac-point} satisfies the following pointwise estimate : 
\begin{align}
r \,  | (\nabla ^{a_0} _{e_i} \alpha _0)&\otimes e^i .\beta_1 |_x^2 \le  \nonumber \\
 \le & Cr\exp \left( - \frac{\sqrr}{C} \text{dist}(x,\alpha _0 ^{-1} (0)) \right) \cdot \exp \left( -
\frac{\sqrr}{C} \text{dist}(x,\alpha _1 ^{-1} (0)) \right)     \label{eq:estimm1}
\end{align}
The second term in \eqref{eq:dirac-point} satisfies the same bound. 
In order for the right hand side of \eqref{eq:estimm1} to pointwise converge to zero, it is 
sufficient and necessary that there exist some $r_0 \ge 1$ such that for all $r\ge r_0$, the
distance from $\alpha _0 ^{-1} (0)$ to  $\alpha _1 ^{-1} (0)$ 
be bounded from below by some $r$-independent $M > 0$.  This condition, under the map
$\Theta$ from \eqref{eq:grsw}, is the Seiberg-Witten equivalent 
of the condition that $\Sigma _i = \Theta (A_i,\psi _i)$ be disjoint curves. Thus, from now
onward we will make the following assumption.

\begin{assumption}  \label{assump:mainassump}  
As above, let $E_0, E_1 \in H^2(X;\mathbb{Z})$ be two line bundles over $X$. Let $(\psi _i, a_i)$, $i=0,1$,  be two solutions to the Seiberg-Witten 
equations \eqref{eq:swgen} for the \spins $W_{E_i}$ with $\psi _i = \sqrt{r} \, (\alpha _i \otimes u_0 ,  \beta _0)$ and $\alpha _i \in \Gamma (E_i)$. We henceforth 
make the assumption that there exists an $r_0\ge 1$ and $M >0$ such that for all $r\ge r_0$ the  inequality  
\begin{equation}
\mbox{dist}(\alpha_0^{-1} (0),\alpha_1^{-1} (0)) \ge  M \label{eq:mainassump} 
\end{equation}
holds. 
\end{assumption}

We now proceed by looking at the second equation in \eqref{eq:swgen}: 
\begin{align}
F& _a ^{1,1} - \frac{i}{8} r \,(|\alpha|^2 -1 - |\beta|^2) \, \omega = \nonumber \\
& = F_{a_0} ^{1,1} + F_{a_1} ^{1,1} - \frac{i}{8} r \, (|\alpha _0|^2 \cdot |\alpha _1| ^2 - 1 - 
|\alpha _0|^2 \cdot |\beta _1|^2 - |\alpha _1| ^2 \cdot |\beta _0| ^2 - \nonumber \\
& \quad \quad \quad \quad - 2 \langle \alpha_0 \beta_1 ,
\alpha _1 \beta_0 \rangle) \, \omega \nonumber \\
& = F_{a_0} ^{1,1} + F_{a_1} ^{1,1} - \frac{i}{8}r\,|\alpha _1|^2 (|\alpha_0|^2 -1 -
|\beta_0|^2)\, \omega  - \frac{i}{8}r\,|\alpha _0|^2 ( |\alpha_1|^2 -1 - 
|\beta_1|^2 )\, \omega  \nonumber \\
& \quad \quad \quad \quad    +  \frac{i}{8}r\, (|\alpha _0|^2-1)(|\alpha _1|^2-1) \, \omega  + 
\frac{i}{4}r\, 
\langle \alpha_0 \beta_1 , \alpha _1 \beta_0 \rangle\, \omega \nonumber \\
& = \frac{i}{8}r\,(1-|\alpha _1|^2) (|\alpha_0|^2 -1 - |\beta_0|^2)\, \omega  -
\frac{i}{8}r\,(1-|\alpha _0|^2) ( |\alpha_1|^2 -1 - |\beta_1|^2 )\, \omega   + \nonumber \\
& \phantom{ = } + \frac{i}{8}r\, (|\alpha _0|^2-1)(|\alpha _1|^2-1)\, \omega   + \frac{i}{4}r\, 
\langle \alpha_0 \beta_1 , \alpha _1 \beta_0 \rangle \, \omega  \nonumber 
\end{align}
From this last equation, and again using \eqref{eq:mainbounds}, one easily deduces that:
\begin{align}
|F_a ^{1,1} -& \frac{i}{8} r \,(|\alpha|^2 -1 - |\beta|^2) \, \omega  | \le  \label{eq:estim2}  \\
&\le Cr\exp \left( - \frac{\sqrr}{C} \text{dist}(x,\alpha _0 ^{-1} (0)) \right) 
\cdot \exp \left( - \frac{\sqrr}{C} \text{dist}(x,\alpha _1 ^{-1} (0)) \right)  + \frac{C}{\sqrr}
\nonumber
\end{align}
Finally, we consider the third equation in \eqref{eq:swgen}: 
\begin{align}
F_a^{0,2}& - \frac{i}{4} r\, \bar{\alpha} \beta = F_{a_0} ^{0,2} + F_{a_1} ^{0,2} -
\frac{i}{4} r\, \overline{\alpha_0 \alpha _1} 
(\alpha _0 \beta _1 + \alpha _1 \beta _0)  \nonumber  \\
 & = \frac{i}{4}r\, \bar{\alpha _0} \beta _0 + \frac{i}{4}r\, \bar{\alpha _1} \beta _1 -
\frac{i}{4} r \, |\alpha _0|^2 \bar{\alpha _1}\beta _1 -
\frac{i}{4} r \, |\alpha _1|^2 \bar{\alpha _0}\beta _0 \nonumber \\
 & = \frac{i}{4}r\, ( 1-|\alpha _1|^2) \bar{\alpha _0}\beta _0 +  \frac{i}{4}r\, ( 1-|\alpha _0|^2)
\bar{\alpha _1}\beta _1 \nonumber
\end{align}
Once again using the bounds \eqref{eq:mainbounds}, we find from this last equation: 
\begin{align}
|F_a^{0,2}& - \frac{i}{4} r\, \bar{\alpha} \beta| \le   \\ 
 & \le Cr\exp \left( - \frac{\sqrr}{C} \text{dist}(x,\alpha _0 ^{-1} (0)) \right) 
\cdot \exp \left( - \frac{\sqrr}{C} \text{dist}(x,\alpha _1 ^{-1} (0)) \right)  + \frac{C}{\sqrr}
\nonumber
\end{align}
To summarize, we have proved the following result.

\begin{proposition} Let $(a,\psi )$ be defined as in \eqref{eq:def1} and assume that there exists an
$r_0\ge 1$ and $M>0$ such that for all $r\ge r_0$, the distance $dist (\alpha _0 ^{-1} (0), \alpha _1 ^{-1} (0))$ is bounded from below by $M$. 
Then for large enough $r$ and any $x\in X$ the pointwise
bound below holds:
\begin{equation} 
|(D_a (\psi), F_a ^{1,1}-\frac{i}{8}r\, (|\alpha|^2-1-|\beta|^2)\, \omega  , F_a ^{0,2} - \frac{i}{4}r
\bar{\alpha} \beta)|_x \le \frac{C}{\sqrr} \label{eq:impos}
\end{equation}
\end{proposition}
\subsection {Inverting the linearized operators of $(a_i , \psi _i)$} \label{sec:32}

This section serves as a digression. The main result of the section is theorem
\ref{theorem:estimm}, an asymptotic (as $r\rightarrow \infty$) regularity statement for the linear operators $L_{(a_i, \psi _i )}$. 

We start with
two easy auxiliary lemmas: 

\begin{lemma}  \label{lemma:surjective} Let $L:V \rightarrow W$ be a surjective Fredholm
operator between Hilbert spaces. Then there exists a $\delta > 0$ such that 
for every linear operator $\ell :V \rightarrow W$ with $||\ell(x)||_W \le \delta \, ||x||_V$, the
operator $L+\ell$ is still surjective.
\end{lemma}

\begin{proof} Since $L$ is Fredholm, we can orthogonally decompose $V$ as
$V=\mbox{Ker} (L)$ $\oplus\, \mbox{Im} (L^*)$. Let $L_1$ be the 
restriction of $L$ to $\mbox{Im} (L^*)$. Then $L_1:\mbox{Im} (L^*) \rightarrow W$ is an
isomorphism with bounded inverse $L_1^{-1}$. 

If the lemma were not true then we could find for all integers $n\ge 1$ an operator $\ell _n : V
\rightarrow W$ with $||\ell _n x ||_W \le 1/n \cdot ||x||_V$ and 
with Coker$(L+\ell _n )\ne \{ 0 \} $. Let $0\ne y_n \in \mbox{Coker}(L+\ell _n)$ with $||y_n||_W
= 1$ and $x_n = L_1^{-1} (y_n)$. Notice that the sequence 
$\{ x_n \} _n $ is bounded by $||L_1 ^{-1} ||$. Since $y_n \in \mbox{Coker}(L+\ell _n)$, $y_n$
is  orthogonal to Im$(L+\ell _n)$. In particular, 
$$ \langle (L+\ell _n ) x_n , y_n \rangle = 0 $$
This immediately leads to a contradiction for large enough $n$ 
since $\langle L\, x_n , y_n \rangle = 1$ and $|\langle \ell _n x_n , y_n \rangle | \le ||L_1^{-1}||/n$. 
\end{proof} 
\begin{lemma} \label{lemma:aux2} Let $V$ and $W$ be two finite rank vector bundles over $X$ and 
$L_r:L^{1,2}(V)\rightarrow L^2 (W)$ a smooth one-parameter family (indexed by $r\ge 1$) of
elliptic, first order, differential operators of index zero. 
Assume further that there exists a $\delta >0$ and $r_0\ge 1$ such that for any zeroth order
linear operator $\ell : L^{1,2} (V) \rightarrow L^2 (W)$ with 
$|| \ell (x) ||_2 < \delta ||x||_{1,2}$, the operator $L_r + \ell$ is onto. Then there exists a $r_1\ge
r_0$ and a $M >0$ such that for all $r\ge r_1$ the inverses of the operators $L_r$ are 
uniformly bounded by $M$, i.e.\ $||L_r ^{-1} y ||_{1,2} \le M ||y || _2 $. 
\end{lemma}

\begin{proof} Notice that a universal upper bound on $L_r ^{-1}$ is equivalent to a
universal lower bound on $L_r$. 
Suppose the lemma were not true: then there would be a sequence $r_n\rightarrow \infty$ and $x_n\in
L^{1,2}(V)$ with $||x_n||_{1,2} = 1$ and 
$||L_{r_n} \, x_n ||_2 < 1/n$. Choose $n$ large enough so that $1/n < \delta$ and define the
operator $\ell : L^{1,2} (V) \rightarrow L^2 (W)$ by 
$\ell (x) = - \langle x_n , x\rangle_{1,2} \cdot L_{r_n} (x_n) $. For this $\ell$ the assumption of
the lemma is met, namely
$$||\ell (x) ||_2 \le \frac{1}{n} ||x||_{1,2} < \delta \, ||x||_{1,2} $$
Thus the operator $L_{r_n} + \ell$ should be onto and injective (since the index of $L_r+\ell$ is
zero). But $x_n$ is clearly a nonzero kernel element. 
This is a contradiction. 
\end{proof} 

Recall that the set ${\cal J}$ of  almost-complex structures compatible with the symplectic form
$\omega$, contains a Baire subset ${\cal J}_0$ of generic 
almost-complex structures in the sense of Gromov-Witten theory (see \cite{kn:taub4}). Also, as
in the introduction, let 
\begin{equation} \Theta : {\cal M}^{SW}_X (W_E) \rightarrow {\cal M}^{Gr}_X (E) 
\end{equation}
be the map introduced in \cite{kn:taub1} which associates an embedded $J$-holomorphic curve
to a Seiberg-Witten monopole. 
\begin{proposition} \label{prop:imp} Let $J$ be chosen from ${\cal J}_0$ and let $(a,\psi )$ be a
solution of the Seiberg-Witten equations \eqref{eq:swgen} 
such that $\Theta (a,\psi) $ doesn't contain any multiply covered components. 
Then there exists a $\delta > 0$ and an $r_0 \ge 1$ such 
that for all linear operators $\ell : L^{1,2} (i \Lambda ^1 \oplus E\otimes W^+ _0) \rightarrow
L^2 (i\Lambda ^0 \oplus i\Lambda ^{2,+} \oplus 
E\otimes W^- _0 )$ with norm $||\ell (x) ||_2 < \delta ||x||_{1,2}$, the operator $L_{(a,\psi)} + \ell
$ is surjective. 
\end{proposition}
Before proceeding to the proof, notice that  proposition \ref{prop:imp} and lemma  \ref{lemma:aux2} immediately imply 
the following theorem, the main result of this section: 
\begin{theorem} \label{theorem:estimm} Choose $J\in {\cal J}_0$ and let  $(a,\psi)$ be a
solution of the Seiberg-Witten equations for the \spin $W_E$ with parameter $r$.
Assume that $\Theta (a,\psi )$ contains no multiply covered components.  
Then there exists a $r-$independent  $M > 0$ and $r_0\ge 1$ such that for all $r\ge r_0$:
\begin{equation}
|| L^{-1} _{(a,\psi)} x ||_{1,2} \le M ||x||_2 \label{eq:contr2}   
\end{equation}
\end{theorem}
\begin{proof}[Proof of proposition \ref{prop:imp}]  The proof is a bit technical and relies on the even more technical account from
\cite{kn:taub2} on the connection between the 
deformation theory of the Seiberg-Witten equations on one hand and the Gromov-Witten
equation on the other. The idea is however very simple:
for large $r\gg 1$, a certain perturbation of the operator $L$ (with the size of the perturbation
getting smaller with larger $r$) has no cokernel if a certain 
perturbation of the linearization of the generalized del-bar operator has no cokernel. The latter is
ensured by the choice of a generic almost complex 
structure $J$ from the Baire set ${\cal J}_0$ of almost complex structures compatible with
$\omega$.

Before proceeding, the (interested) reader is advised to familiarize him/her-self with the notation from 
\cite{kn:taub2}, in particular, sections 4 and 6 as the remainder of the proof heavily relies on it. 
For convenience we restate here the parts of lemma 4.11 and a slightly
modified version of lemma 6.7 from \cite{kn:taub2} relevant to our situation.

\medskip
{\bf Lemma 4.11}\qua {\sl The equation $Lq+\eta q = g$ is solvable if and only if, for each
$k$ 
$$\Delta_{c_k} w^k + \gamma _0 ^k (w) + \eta _k (w) = x (g^k) + \gamma _1^k (g) $$
}

\medskip
Notice that the assignment of $\eta _k$ to $\eta $ is linear i.e.\ for two operators $\eta$ and $\eta
'$, we have 
$(\eta + \eta ')_k = \eta _k + \eta ' _k$.  

\medskip
{\bf Lemma 6.7$\,'$}\qua {\sl  The equation $(L_{\Psi _r (y)} + \ell ) p = g$ has an
$L^{1,2}$ solution $p$ if and only if there exists 
$u=(u^1,...,u^k)\in \oplus _k L^{1,2} (N^{(k)})$ for which 
$$  \Delta_y u^k + \phi _0 ^k (u) + \ell _k (u) = \Upsilon _1^{-1} x (g^k) + \phi _1^k (g) $$
holds for each $k$. }
 \medskip

The proof of lemma 6.7$\,'$ is almost identical to that of the original lemma 6.7 in
\cite{kn:taub2}. The only difference is in {\sl Step 2} where 
Taubes shows that one can write the equation $L_{\Psi _r (y)} p = g$ in the form $Lp + \eta p = 
g$ with $L$ as in lemma 4.11 above and with $\eta$ an appropriate (bounded) correction term
(see (6.30) in \cite{kn:taub2} for a precise definition). 
The difference here is that in our case one can write $(L_{\Psi _r (y)} + \ell ) p = g$ as $Lp +
\eta ' p =  g$ (with $L$ again as in lemma 4.11 of \cite{kn:taub2}) but 
with $\eta '  = \eta + \ell$. Since $\ell$ is assumed bounded, lemma 4.11
applies to $\eta '$ in the exact same  way as it applied 
to the original $\eta$ and the proof of lemma 6.7 in \cite{kn:taub2} transfers verbatim to our
case. Note also that the operators $\phi _i ^k$ occurring 
in lemmas 6.7 and 6.7$\,'$ are identical so in particular they continue to satisfy the bounds
asserted by lemma 6.7 of \cite{kn:taub2}. 

According to  lemma \ref{lemma:surjective} there exists a $\delta ' > 0$ such that $\Delta _y +
\ell '$ is still surjective if $||\ell '|| < \delta '$. Choose $r$ large 
enough so that  $||\phi _0 ^k||<\delta '/2k$. On the other hand, since $\ell _k (v) = \pi (\chi
_{25\delta,k} \ell (\sum _{k'} \chi _{100\delta,k'} \underline{v}^{k'}))  $
we find that $||\ell _k|| \le C\, ||\ell ||$. Thus choosing $\delta = \delta ' /2C$ ensures that $ L_{\Psi
_r (y)} + \ell$ is surjective provided that $||\ell||<\delta$. 
This finishes the proof of proposition \ref{prop:imp}. 
\end{proof}
\subsection {The linearized operator at $(a,\psi)$} \label{section:four}
In order to use the contraction mapping principle to deform the approximate solution $(a,\psi )$
to an honest solution of the 
Seiberg-Witten equations, we need to know that $L=L_{(a,\psi)}$ admits an inverse whose norm
is bounded independently 
of $r$.  We start by  exploring when the equation 
\begin{equation}
L\xi = g \label{eq:main}
\end{equation}
has a solution $\xi$ for a given $g$. Here: 
$$\xi \in L^{1,2} (i\Lambda ^1 \oplus (E_0\otimes E_1\otimes W^+ _0))\quad \mbox{and}
\quad g \in 
L^2 (i\Lambda ^0 \oplus i\Lambda ^{2,+} \oplus (E_0\otimes E_1\otimes W^- _0))$$
The idea is to restrict equation \eqref{eq:main} first to a neighborhood of $\alpha _0 ^{-1}(0)$.
Over such a neighborhood the bundle $E_1$ is trivial and,
under an isomorphism trivializing $E_1$, the  
equation \eqref{eq:main} becomes a zero-th order perturbation of the equation $L_0 \xi_0  =g_0
$ (with $\xi_0$ and $g_0$ being appropriately defined in terms of $\xi$ and $g$). This allows one  
to take advantage of the results of theorem \ref{theorem:estimm} about the inverse of $L_0=L_{(a_0, \psi _0)}$.
Then one restricts \eqref{eq:main} to a neighborhood of $\alpha _1 ^{-1}(0)$ where 
the bundle $E_0$ trivializes and once again uses theorem \ref{theorem:estimm}, this time for
the inverse of $L_1=L_{(a_1, \psi_1)}$. Finally, one restricts to the 
complement of a neighborhood of 
$\alpha _0 ^{-1}(0) \cup \alpha _1 ^{-1}(0)$ where both $E_0$ and $E_1$ become trivial and
$L$ becomes close to $S$ - the linearized operator of the unique 
solution $({\cal A}_0, \sqrr \, u_0)$ for the anticanonical \spin $W_0$. 

To begin this process, choose regular neighborhoods $V_i$ of $\alpha _i ^{-1} (0)$, $i=0$, $1$ 
subject to the condition:
$$\mbox{dist} (V_0 , V_1) \ge M \quad \quad \mbox{for some } M>0 $$ 
The existence of such neighborhoods $V_i$ follows from our main assumption
\eqref{eq:mainassump}. A priori, as one chooses larger values of $r$, it seems
that the sets $V_i$ may need to be chosen anew as well. However, it was shown in
\cite{kn:taub1}, section 5c, that in fact this is not necessary. An initial ``smart" choice
of $V_i$ for large enough $r$ ensures that for $r'>r$, the zero sets $\alpha _i ^{-1} (0)$ continue
to lie inside of $V_i$.
Choose an open set $U$ such that $X=V_0\cup V_1\cup U$ and such that 
$$  U \cap (\alpha _0 ^{-1}(0) \cup \alpha _1 ^{-1}(0)) = \emptyset $$
Arrange the choices of $V_i$ and $U$ further so that $\partial V_i$ is an embedded 3-manifold
of $X$ and so 
that $U\cap V_i $ contains a collar $\partial V_i \times I$. Here $I$ is some segment $[0,d]$ and
$\partial V_i$ corresponds to $\partial V_i \times \{ d\}$.  For the sake of simplicity of notation, we shall make the
assumption that for large values of $r$, the sets $\alpha ^{-1} _i (0)$, $i=0,1$, are connected. The case of disconnected zero sets of the 
$\alpha _i$'s is treated much in the same way except for that in the following, one would have to choose a bump function $\chi _{\delta, i}$ 
(see below) for each connected component. This complicates notation to a certain degree but doesn't lead to new phenomena.  

Fix once and for all a bump function $\chi : [0,\infty ) \rightarrow [0,1]$ which is 1 on [0,1] and
0 on $[2,\infty )$. 
For $0<\delta <d/1000 $ define 
$\chi _{\delta , i } : X  \rightarrow [0,1]$ by: 
\begin{equation}
\chi _{\delta ,i } (x) = \left\{ \begin{array}{ll}
                                 1 & x\in V_i \backslash (\partial V_i \times I)  \\
                                 \chi (t/ \delta) & x=(y,t)\in \partial V_i \times I \\
                                  0 & x\ \not\in V_i 
\end{array}
\right.
\end{equation}
Set $V_0'=V_0\cup U$ and $V_1' = V_1 \cup U$. Define the isomorphisms $\Upsilon _0 :
\mathbb{C}\times V_0 ' \rightarrow E_1\big| _{V'_0}$ and 
$\Upsilon _1 : \mathbb{C}\times V_1 ' \rightarrow E_0\big| _{V'_1}$
 as 
$\Upsilon _0 (\lambda , x) = \alpha _1 (x) \cdot \lambda$ and 
$\Upsilon _1 (\lambda , x) = \alpha _0 (x) \cdot \lambda$.  
Also, for $i=0$, $1$ define the operators 
\begin{align*}
M_i : L^{1,2} (i\Lambda ^1 \oplus (E_i\otimes W^+ _0);V_i ')& \rightarrow  L^2 (i\Lambda ^0
\oplus i\Lambda ^{2,+} \oplus (E_i\otimes W^- _0);V_i ') \\
& \text{ and }  \\ 
T: L^{1,2} (i\Lambda ^1 \oplus W^+ _0;U) & \rightarrow  L^2 (i\Lambda ^0 \oplus i\Lambda
^{2,+} \oplus W^- _0;U)
\end{align*}
 by demanding the diagrams
$$\begin{CD}
 L^{1,2} (i\Lambda ^1 \oplus (E_0\otimes E_1 \otimes W^+ _0);V_i ')   @<\Upsilon _i<< 
L^{1,2} (i\Lambda ^1 \oplus W^+ _i;V_i ') \\
@V{L} VV     @VV{M_i}V  \\
L^2 (i\Lambda ^0 \oplus i\Lambda ^{2,+} \oplus (E_0\otimes E_1\otimes W^- _0);V_i ')
@<\Upsilon _i<<  
L^2 (i\Lambda ^0 \oplus i\Lambda ^{2,+} \oplus W^- _i;V_i ') 
\end{CD}$$and $$\begin{CD}
 L^{1,2} (i\Lambda ^1 \oplus (E_0\otimes E_1 \otimes W^+ _0);U)   @<{\Upsilon _0 \circ
\Upsilon _1}<<  L^{1,2} (i\Lambda ^1 \oplus W^+ _0;U) \\
@V{L} VV     @VV{T}V  \\
L^2 (i\Lambda ^0 \oplus i\Lambda ^{2,+} \oplus (E_0\otimes E_1\otimes W^- _0);U)
@<{\Upsilon _0 \circ \Upsilon _1}<<  
L^2 (i\Lambda ^0 \oplus i\Lambda ^{2,+} \oplus W^- _0 ;U) 
\end{CD}$$
to be commutative diagrams (with $W^\pm _i = E_i \otimes W^\pm _0$).  

We now start our search for a solution $\xi$ of \eqref{eq:main}  in the form: 
\begin{equation}
 \xi = \Upsilon _0 (\chi _{100\delta,0}\xi_0) + \Upsilon _1 (\chi _{100\delta,1}\xi _1) + \Upsilon
_0 \Upsilon _1 \left( 
(1-\chi _{4\delta ,0})(1-\chi _{4\delta, 1}) \eta \right)   \label{eq:form}
\end{equation}
Here $\xi _i \in  L^{1,2} (i\Lambda ^1 \oplus (E_i \otimes W^+ _0))$ and $\eta \in L^{1,2}
(i\Lambda ^1 \oplus W^+ _0) $. Given a 
$g\in L^2 (i\Lambda ^0 \oplus i\Lambda ^{2,+} \oplus (E_0\otimes E_1\otimes W^- _0))$,
define 
$g_i \in L^2 (i\Lambda ^0 \oplus i\Lambda ^{2,+} \oplus (E_i\otimes W^- _0))$ and $\gamma
\in L^2 (i\Lambda ^0 \oplus i\Lambda ^{2,+} \oplus W^- _0)$ as 
\begin{equation}
g_i = \Upsilon _i ^{-1} (\chi _{25\delta,i} g) \quad  \text{and} \quad \gamma = (\Upsilon _0
\Upsilon _1)^{-1} 
\left( (1-\chi _{25\delta ,0})(1-\chi _{25\delta, 1}) g \right) \label{eq:formm}
\end{equation} 
It is easy to check that $g$, $g_i$ and $\gamma$ satisfy a relation similar to \eqref{eq:form},
namely: 
\begin{equation}
 g = \Upsilon _0 (\chi _{100\delta,0}g_0) + \Upsilon _1 (\chi _{100\delta,1}g _1) + \Upsilon _0
\Upsilon _1 \left( 
(1-\chi _{4\delta ,0})(1-\chi _{4\delta, 1}) \gamma \right)   \label{eq:formmm}
\end{equation}
Putting the form \eqref{eq:form}  of $\xi$ and the form \eqref{eq:formmm} of $g$ 
into equation \eqref{eq:main}, after a few simple manipulations, yields the equation
\begin{align}
\Upsilon_0 & (\chi _{100\delta,0}(M_0 (\xi _0) -\Upsilon _1 {\cal P} (d\chi _{4\delta,0} , \eta) -
g_0)) +  \label{eq:split}  \\ 
+ & \Upsilon_1 (\chi _{100\delta,1}(M_1 (\xi _1) -\Upsilon _0 {\cal P} (d\chi _{4\delta,1} ,
\eta) - g_1)) + &  \nonumber  \\
+ & \Upsilon _0 \Upsilon _1 ((1-\chi _{4\delta,0})(1-\chi _{4\delta,1})\left( T\eta + \Upsilon_1 ^{-1} {\cal
P}(d\chi _{100\delta,0},\xi _0) + \right.  \nonumber \\
+ & \left. \Upsilon _0 ^{-1} {\cal P}(d\chi _{100\delta,1},\xi _1) - \gamma \right)  \nonumber \\ 
  = & \, \, 0 \nonumber 
\end{align}
In the above, ${\cal P}$ denotes the principal symbol of $L$. 
This last equation suggests a splitting into three equations (each corresponding to one line in
\eqref{eq:split}):
\begin{align}
M_0(\xi _0) -  \Upsilon _1 {\cal P} (d\chi _{4\delta,0} , \eta) = g_0 \nonumber \\
M_1 (\xi _1) -\Upsilon _0 {\cal P} (d\chi _{4\delta,1} , \eta) =  g_1   \label{eq:splitt}  \\
T\eta + \Upsilon _1 ^{-1}{\cal P}(d\chi _{100\delta,0},\xi _0) +  \Upsilon _0 ^{-1}{\cal P}(d\chi _{100\delta,1},\xi _1) = \gamma
\nonumber 
\end{align}
Equation \eqref{eq:split} (and hence also equation \eqref{eq:main}) can be recovered from
\eqref{eq:splitt}  by multiplying the three equations by 
$\Upsilon_0 \cdot \chi _{100\delta,0}$,  $\Upsilon_1 \cdot \chi _{100\delta,1}$ and $\Upsilon
_0 \Upsilon _1 \cdot ((1-\chi _{4\delta,0})(1-\chi _{4\delta,1})$ 
respectively and then adding them. Thus, given a $g$ and with $g_i$ and $\gamma$ defined by
\eqref{eq:formm},  solutions $\xi _i$ and 
$\eta$ of \eqref{eq:splitt} lead to a solution $\xi$ of \eqref{eq:main} via \eqref{eq:form}. 
However, the problem with \eqref{eq:splitt} is that the operators $M_i$ and $T$ are not
defined over all of $X$. We remedy this in the next step. 

Define new operators:
\begin{align*}
M_i' : L^{1,2} (i\Lambda ^1 \oplus (E_i\otimes W^+ _0))& \rightarrow  L^2 (i\Lambda ^0
\oplus \Lambda ^{2,+} \oplus (E_i\otimes W^- _0)) \\
& \text{ and }  \\ 
T': L^{1,2} (i\Lambda ^1 \oplus W^+ _0) & \rightarrow  L^2 (i\Lambda ^0 \oplus \Lambda
^{2,+} \oplus W^- _0)
\end{align*}
 by 
\begin{align}
M_i'= \chi _{200\delta,i} M_i + (1-\chi _{200\delta,i}) L_i \quad \quad \nonumber  \\
 T' = (1-\chi _{\delta,0})(1-\chi _{\delta,1})T + (\chi _{\delta,0} + \chi _{\delta,1}) S
\label{eq:formmmm}
\end{align}
Here $L_i=L_{(a_i,\psi _i)}$. 
Now replace the coupled equations \eqref{eq:splitt} by the following system:
\begin{align}
M_0'(\xi _0) -  \Upsilon _1 {\cal P} (d\chi _{4\delta,0} , \eta) = g_0 \nonumber \\
M_1' (\xi _1) -\Upsilon _0 {\cal P} (d\chi _{4\delta,1} , \eta) =  g_1   \label{eq:splitttt}  \\
T'\eta + \Upsilon _1 ^{-1} {\cal P}(d\chi _{100\delta,0},\xi _0) +  \Upsilon _0 ^{-1} {\cal P}(d\chi _{100\delta,1},\xi _1) = \gamma
\nonumber 
\end{align}
The advantage of \eqref{eq:splitttt} over \eqref{eq:splitt}  is that the former is defined over all
of $X$ (notice that the support of ${\cal P}(d\chi _{100\delta,0},\xi _0)$ lies in the domain of $\Upsilon_1 ^{-1}$ and the 
support of ${\cal P}(d\chi _{100\delta,1},\xi _1)$ lies in the domain of $\Upsilon _0 ^{-1}$). On the other hand, solutions of \eqref{eq:splitttt} give 
rise to solutions of \eqref{eq:main}  in the same way as solutions of \eqref{eq:splitt} did because 
\begin{align}
\chi _{100\delta,i} \cdot M_i '= \chi _{100\delta,i} \cdot M_i  \quad \quad i=0,\, 1 \nonumber \\
(1-\chi _{4\delta,0}) (1- \chi _{4\delta ,1}) T' = (1-\chi _{4\delta,0}) (1- \chi _{4\delta ,1}) T
\nonumber 
\end{align}
\begin{lemma} \label{lemma:goodone} For every $\epsilon>0$ there exists an $r_\epsilon \ge 1$ such that for $r\ge
r_\epsilon$ the following hold:
\begin{align}
|| (M_i ' - L_i) x_i   \,||_2 & \le \epsilon ||x _i ||_2  \nonumber  \\  
||(T'-S)y \, ||_2   & \le \epsilon ||y||_2     \nonumber
\end{align} 
Here $x_i \in L^{1,2} (i\Lambda ^1 \oplus E_i \otimes W^+_0 )$ and $y\in  L^{1,2} (i\Lambda
^1 \oplus W^+_0 )$.
\end{lemma}
\begin{proof} The above Sobolev inequalities are proved by first calculating pointwise bounds for
$| (M_i ' - L_i) x_i   \,|_p$ and $|(T'-S)y \, |_p$, $p\in X$.
Notice firstly that $|(M_i ' - L_i) x_i   \,|_p=0$ if $p\notin V_i$ and $|(T'-S)y \, |_p=0$ if $p \notin
U$. For $p\in V_i$ and for $q\in U$, a straightforward 
but somewhat tedious calculation
shows that: 
\begin{align}
| (M_i ' - L_i) x_i   \,|_p & \le C\, \left( \sqrr\, |1-|\alpha _i |^2| + \sqrr \, |\beta _i|\, |\alpha _i| +
|\nabla ^{a_i} \alpha _i|\right) |x_i|_p    \nonumber   \\ 
|(T'-S)y \, |_q & \le C\, (   \sqrr \,  |1-|\alpha _0|^2| + \sqrr \,  |1-|\alpha _1|^2| + \sqrr \, |\beta _0| +
\nonumber \\ 
       & \phantom{mmmm} + \sqrr \, |\beta _1| +|\nabla ^{a_0} \alpha _0| +
 |\nabla ^{a_1} \alpha _1|  ) \, |y|_q \nonumber   
\end{align}
Squaring and then integrating both sides over $X$ together with a reference to
\eqref{eq:mainbounds} gives the desired Sobolev inequalities.
\end{proof}

The lemma suggests that the system \eqref{eq:splitttt} can be replaced by  the system:
\begin{align}
L_0(\xi _0') -  \Upsilon _1 {\cal P} (d\chi _{4\delta,0} , \eta') = g_0 \nonumber \\
L_1 (\xi _1') -\Upsilon _0 {\cal P} (d\chi _{4\delta,1} , \eta') =  g_1   \label{eq:splittttt}  \\
S\eta' + \Upsilon _1 ^{-1} {\cal P}(d\chi _{100\delta,0},\xi _0') +  \Upsilon _0 ^{-1} {\cal P}(d\chi _{100\delta,1},\xi _1') = \gamma
\nonumber 
\end{align}
Lemmas \ref{lemma:goodone} and \ref{lemma:surjective} say that for $r\gg 0$, \eqref{eq:splitttt} has a solution $(\xi _0, \xi _1, \eta)$  if 
 \eqref{eq:splittttt} has a solution $(\xi _0', \xi _1', \eta')$. It is this latter set of equations that we now 
proceed to solve. 

Since $S$ is onto, we can solve the third equation in \eqref{eq:splittttt}, regarding $\xi _0'$ and
$\xi _1'$
as parameters. Thus: 
\begin{equation}
\eta' = \eta' (\xi _0', \xi _1') = S^{-1} (\gamma - \Upsilon _1 ^{-1} {\cal P}(d\chi _{100\delta,0},\xi _0') -  \Upsilon _0 ^{-1}{\cal
P}(d\chi _{100\delta,1},\xi _1') ) \label{eq:estimm}
\end{equation}
Recall that the inverse of $S$ satisfies the bound \eqref{eq:estim}:
\begin{equation} 
||S^{-1} y || _2 \le \frac{4}{\sqrr} ||y||_2 \quad \text{for} \quad 
y\in L^2 (i\Lambda ^0 \oplus i\Lambda ^{2,+} \oplus W^- _0) \nonumber
\end{equation}
We will solve the first two equations in \eqref{eq:splittttt} simultaneously by first rewriting
them in the form:
\begin{align}
\xi _ 0' = & L_0^{-1} (g_0 +  \Upsilon _1  {\cal P}(d\chi _{4\delta,0}, \eta'(\xi _0',\xi _1')))
\nonumber \\
\xi _ 1' = & L_1^{-1} (g_1 + \Upsilon _0 {\cal P}(d\chi _{4\delta,1}, \eta'(\xi _0',\xi _1')))
\label{eq:contr}
\end{align} 
To solve \eqref{eq:contr}  is the same as to find a fixed point of the map 
$Y: L^{2} (i\Lambda ^1\oplus W^+ _{E_0}) \times L^2 (i\Lambda ^1\oplus W^+ _{E_1})   
 \rightarrow L^2  (i\Lambda ^1\oplus W^+ _{E_0}) \times L^2 (i\Lambda ^1\oplus W^+
_{E_1})  $ given by
\begin{align}
Y(\xi _0',\xi _1') &  = \label{eq:llast} \\
= &  (L_0^{-1} (g_0 + \Upsilon _1  {\cal P}(d\chi _{4\delta,0},
\eta')) , 
L_1^{-1} (g_1 + \Upsilon _0 {\cal P}(d\chi _{4\delta,1}, \eta') ))  \nonumber
\end{align}
with $\eta '$ given by \eqref{eq:estimm}. 
The existence and uniqueness of such a fixed point will be guaranteed by the fixed point
theorem for Banach spaces if we can show that $Y$ is a contraction mapping. 
To see this, let $x$, $y\in L^2  (i\Lambda ^1\oplus W^+ _{E_0}) \times L^2 (i\Lambda ^1\oplus
W^+ _{E_1}) $ be two arbitrary sections. Using the first bound of \eqref{eq:estim} and the
result of theorem  \ref{theorem:estimm} 
to bound the norms of $L_i^{-1}$, one finds: 
\begin{align}
||Y&(x)-Y(y)||^2 _{2}  = \nonumber \\
=\, & ||L_0^{-1} (g_0 + \Upsilon _1  {\cal P}(d\chi _{4\delta,0}, \eta(x))) - 
L_0^{-1} (g_0 + \Upsilon _1 {\cal P}(d\chi _{4\delta,0}, \eta(y)))||^2 _{2} \nonumber \\
+\, & ||L_1^{-1} (g_1 + \Upsilon _0 {\cal P}(d\chi _{4\delta,1}, \eta(x)) ) - L_1^{-1} (g_1 +
\Upsilon _0 {\cal P}(d\chi _{4\delta,1}, \eta(y) ))||^2 _{2} \nonumber \\
\le & C_0 ||\eta (x) - \eta (y) ||^2 _{2} + C_1 ||\eta (x) - \eta (y) ||^2 _{2} \nonumber \\
\le & C ||S^{-1} (\Upsilon _1 ^{-1} {\cal P}(d\chi _{100\delta,0},y) - \Upsilon _1 ^{-1} {\cal P}(d\chi _{100\delta,0},x) + \\ 
& + \Upsilon _0 ^{-1}{\cal P}(d\chi _{100\delta,1},y) - \Upsilon _0 ^{-1} {\cal P}(d\chi _{100\delta,1},x)) ||^2 _{2} \nonumber \\
\le & \frac{C}{r} ||x-y||^2_2  \label{eq:last} 
\end{align}
Choosing $r>2C$, where $C$ is the constant in the last line of \eqref{eq:last}, makes $Y$ a
contraction mapping. Thus we finally arrive at an $L^2$ solution 
$(\xi_0' , \xi _1')$. It is in fact an $L^{1,2}$ solution because of \eqref{eq:contr}. This, together
with equation \eqref{eq:estimm} provides a solution $(\xi _0', \xi _1', \eta ')$ of \eqref{eq:splittttt}. 
As explained above, this gives rise to a solution $(\xi _0, \xi _1, \eta)$ of \eqref{eq:splitttt} and thus provides a solution 
$\xi \in L^{1,2}  \in(i\Lambda ^1\oplus W^+ _0) $ of \eqref{eq:main}. 
 In particular, we have proved half of the following theorem.
\begin{theorem} \label{theorem:mainn} Let $(a, \psi)$ be constructed from $(a_i ,\psi _i)$ as in
\eqref{eq:def1}. Suppose that the $(a_i, \psi _i)$ meet assumption \ref{assump:mainassump},
that $\Theta(a_i,\psi _i)$ contains no multiply covered tori and that $J$ has been chosen from the Baire set ${\cal J}_0$ of
compatible almost complex structures.  
Then  $L_{(a,\psi )} :L^{1,2} (i\Lambda ^1\oplus E_0\otimes E_1\otimes  W^+ _0)  \rightarrow
L^2(i\Lambda ^0\oplus i\Lambda ^{2,+} \oplus E_0\otimes E_1\otimes  W^- _0)$ 
is invertible with bounded inverse $||L_{(a,\psi )} ^{-1} y ||_{1,2} \le C \, ||y|| _2 $ for all
sufficiently large $r$. Here $C$ is independent of $r$.
\end{theorem}

\begin{proof} It remains to prove the inequality $||L_{(a,\psi )} ^{-1} y ||_{1,2} \le C \, ||y|| _2 $.
Each of the two lines of \eqref{eq:contr}, together with the bound \eqref{eq:contr2} 
on $L_i^{-1}$,  yields:
\begin{equation}
|| \xi _i' || _{1,2} \le C \, (||g_i||_2 + ||\eta' (\xi _0' , \xi _1' )||_2 ) \label{eq:estim1}
\end{equation}
A bound for the second term on the right-hand side of \eqref{eq:estim1} comes from
\eqref{eq:estimm} and the $L^2$ bound in \eqref{eq:estim}:
\begin{equation}
||\eta' (\xi _0' , \xi _1' )||_2 \le \frac{C}{\sqrr}(||\gamma ||_2 + ||\xi _0'||_2 + ||\xi _1' ||_2 ) 
\label{eq:estimm2}
\end{equation}
Adding the two inequalities \eqref{eq:estim1} for $i=0$, $1$ and using \eqref{eq:estimm2}
gives: 
\begin{equation}
(1-\frac{C}{\sqrr}) \, (||\xi _0'||_{1,2} + ||\xi _1' ||_{1,2} ) \le C\, (||g_0||_2 + ||g_1||_2 +
\frac{1}{\sqrr} ||\gamma ||_2) \label{eq:estimm3}
\end{equation}
For large enough $r$, this last inequality gives a bound on the $L^{1,2}$ norm of $(\xi _0' , \xi
_1' )$ in terms of an $r$-independent multiple of the 
$L^2$ norm of $(g_0, g_1, \gamma)$. With this established, the missing piece, namely the
$L^{1,2}$ bound of $\eta'$, comes from 
\eqref{eq:estimm} and the $L^{1,2}$ bound in \eqref{eq:estim}:
\begin{equation}
||\eta'||_{1,2} \le C\, (||\gamma ||_2 + ||\xi _0'||_2 + ||\xi _1' ||_2 ) \le C \, (||\gamma ||_2 + ||g _0||_2 +
||g _1 ||_2 ) \label{eq:estimm4}
\end{equation}
It remains to relate the now established bound on $(\xi _0' , \xi_1', \eta ')$ to a bound for $(\xi _0 , \xi_1, \eta )$. 
To begin doing that, write the systems \eqref{eq:splittttt} and \eqref{eq:splitttt} schematically as: 
$$ {\cal F}(\xi _0' , \xi_1', \eta ') = (g_0,g_1,\gamma) \quad \mbox{ and } \quad   {\cal G}(\xi _0 , \xi_1, \eta ) = (g_0,g_1,\gamma) $$ 
Lemma \ref{lemma:goodone} implies that for any $\varepsilon > 0$ there exists a $r_\varepsilon \ge 1$ such that for all $r\ge r_\varepsilon$ the inequality 
$|| ({\cal F} - {\cal G}) \, x ||_2 \le \varepsilon\, ||x|| _2$ holds. The established surjectivity of ${\cal F}$ guarantees (by means of lemma \ref{lemma:surjective})
that ${\cal G}$ is also surjective. The proof of theorem \ref{theorem:main} thus far, also shows that $||{\cal F}^{-1}||  \le C$ where $C$ is $r-$independent. Now the 
standard inequality 
$$ ||{\cal G}^{-1}|| \le ||{\cal F}^{-1} || + ||{\cal G}^{-1} - {\cal F}^{-1} || \le ||{\cal F}^{-1} || + ||{\cal F}^{-1} ||\cdot ||{\cal G}^{-1} ||\cdot ||{\cal G} - {\cal F} || $$
implies the $r-$independent bound for $||{\cal G}^{-1} ||$: 
$$  ||{\cal G}^{-1}|| \le \frac{||{\cal F}^{-1} ||}{1- ||{\cal F}^{-1} || \cdot ||{\cal G} - {\cal F} || } \le \frac{C}{1-C\varepsilon} $$
This last inequality provides $L^{1,2}$ bounds on $(\xi _0 , \xi _1 )$ and $\eta$ in terms of the $L^2$ norms of $(g_0 ,
g_1)$ and $\gamma$ which in turn imply an $r-$independent  $L^{1,2}$ bound 
on $\xi = L^{-1} g$ in terms of the $L^2$ norm of $g$ through \eqref{eq:form} and
\eqref{eq:formm}. This finishes the proof of theorem \ref{theorem:mainn}. 
\end{proof}
\subsection {Deforming $(a,\psi)$ to an honest solution} \label{section:six} 
The goal of this section is to show that the approximate solution $(a,\psi )$ can be made into an
honest solution of the Seiberg-Witten equations by a  deformation whose size goes to zero as $r$ goes to infinity. 

To set the stage, let $SW:L^{1,2}(i\Lambda ^1 \oplus W_E ^+ ) \rightarrow L^2(i\Lambda ^0
\oplus i\Lambda ^{2,+} \oplus W_E^-)$ 
denote the Seiberg-Witten operator 
$$SW(b, \phi ) = (d^* b\, ,  F_b^+ -F^+_{{\cal A}_0} - q(\phi ,\phi ) + \frac{ir}{8}\omega \, , D_b \phi )$$
We will search for a zero of $SW$ of the form $(a,\psi ) + (a' , \psi ')$ with $(a' , \psi ')\in
B(\delta)$. Here $B(\delta)$ is the closed  ball in 
$L^{1,2}(i\Lambda ^1 \oplus W_E ^+ )$ centered at zero and with radius $\delta >0$ which we
will choose later but which should be thought of as being small. 
The equation $SW((a,\psi ) + (a' ,\psi ')) = 0$ can be written as: 
\begin{equation}
0= SW(a,\psi ) + L_{(a,\psi)} (a' , \psi ') + Q(a' , \psi ') \label{eq:main2}
\end{equation}
Here $Q:L^{1,2} (i\Lambda ^1 \oplus W_E ^+ ) \rightarrow L^2 (i\Lambda ^0 \oplus i\Lambda
^{2,+} \oplus W_E^-) $ is the quadratic map given by: 
\begin{equation}
Q(b,\phi_0,\phi _2) = (b.(\phi_0 + \phi _2) , \frac{i}{8}(|\phi_0|^2 - |\phi _2|^2) \omega +
\frac{i}{4} (\bar{ \phi _0} \, \phi _2 + \phi _0 \, \bar{\phi _2} ))  \label{eq:mapq}
\end{equation}
\begin{lemma} \label{lemma:qlem} 
For $x$, $y\in L^{1,2} (i\Lambda ^1 \oplus W_E ^+ )$, the map $Q$ satisfies the inequality:
\begin{equation}
||Q(x) - Q(y) ||_2 \le C\,  (||x||_{1,2} + ||y||_{1,2} )\, ||x-y||_{1,2} \label{eq:qbound}
\end{equation}
\end{lemma}

\begin{proof} This is a standard inequality for quadratic maps and it can be explicitly checked
using the definition of $Q$ and the 
multiplication theorem for Sobolev spaces. We give the calculation for the 
first component of the right hand side of \eqref{eq:mapq}.  
Let $x=(b,\phi )$ and $y=(c,\varphi)$, then we have:
\begin{align}
||b.\phi - c.\varphi ||_2 = &  ||b.\phi - c.\phi + c.\phi - c.\varphi ||_2 \le ||(b-c).\phi||_2 + ||c.(\phi -
\varphi)||_2 \nonumber  \\
                                    & \le C\, ||b-c||_{1,2}\, ||\phi||_{1,2} + C\,||c||_{1,2} \, ||\phi - \varphi ||_{1,2} 
\nonumber \\
                                    & \le C \, (||(b,\phi) - (c,\varphi) ||_{1,2} ) \, (||(b,\phi)||_{1,2} + ||(c,\varphi
)||_{1,2} ) \nonumber
\end{align}
The other components are checked similarly. 
\end{proof}

Solving equation \eqref{eq:main2} for $(a' , \psi ') \in L^{1,2} (i\Lambda ^1 \oplus W_E ^+ )$ is
equivalent to finding a fixed point of the map $Y:B(\delta ) \rightarrow B(\delta )$ given by:
\begin{equation}
Y(b,\phi ) = - L_{(a,\psi )} ^{-1} (SW(a,\psi ) + Q(b ,\phi)) \label{eq:defofy}
\end{equation}
In order for the image of $Y$ to lie in $B(\delta)$ we need to choose $r$ large enough and
$\delta $ small enough. To make this precise, let 
$(b, \phi )\in B(\delta)$. Using the bounds in \eqref{eq:impos} we find that 
$$||SW(a, \psi ) ||_2 \le \frac{C}{\sqrr} $$
and so together with the results of theorem \ref{theorem:mainn} and lemma \ref{lemma:qlem}
we get:  
\begin{align}
||Y(b, \phi )||_{1,2} \le  \frac{C}{\sqrr} + C\cdot \delta ^2 \nonumber
\end{align}
Choosing $\delta < 1/2C$ and $r> 4C^2/\delta ^2 $ ensures that $Y$ is well defined. 

\begin{lemma} The map $Y:B(\delta) \rightarrow B(\delta)$ as defined by \eqref{eq:defofy} is a
contraction mapping for $r$ large enough and 
$\delta$ small enough. 
\end{lemma}
\begin{proof}  Let $x$, $y\in B(\delta)$, then using \eqref{eq:qbound} we find:
\begin{equation}
||Y(x) - Y(y) ||_{1,2} \le C\, ||Q(x) - Q(y)||_2 \le C \,  ||x+y||_{1,2} \, ||x-y||_{1,2} 
\end{equation}
Choosing $\delta < 1/2C$ makes $C\, ||x+y||_{1,2} \le 2C\delta$ less than 1.  
\end{proof}

We summarize in the following:
\begin{theorem} \label{theorem:main}  Let $(a, \psi)$ be constructed from $(a_i ,\psi _i)$ as in
\eqref{eq:def1}. Suppose that the $(a_i, \psi _i)$ meet assumption \ref{assump:mainassump}, 
that $\Theta(a_i,\psi _i)$ contains no multiply covered tori and that $J$ has been chosen from the Baire set ${\cal J}_0$ of
compatible almost complex structures.  

Then there exists a  $\delta_0 > 0$ such that for any $0< \delta \le \delta _0$ there exists an
$r_{\delta} \ge 1$ such that for every $r\ge r_{\delta}$ 
 there exists a unique solution $(a,\psi ) + (a' , \psi ')$  
of the Seiberg-Witten equations (with perturbation parameter $r$) with $(a', \psi ') \in L^{1,2} (i\Lambda ^1 \oplus W_E^+)$ satisfying the bound 
$||(a' , \psi ')||_{1,2} \le \delta$.  
\end{theorem}
\begin{remark} It is not known if theorem \ref{theorem:main} holds under the relaxed hypothesis allowing $\Theta(a_i,\psi _i)$ to contain multiply 
covered tori. The difficulty in dealing with this case stems from the fact that the operators $L_{(a_i,\psi_i)}$ may no longer have trivial cokernel. 

\end{remark}
\section{Comparison with product formulas} \label{sec:4} 
Before proceeding further, we would like to take a moment to point out the similarities and
differences between our construction of $(A,\psi)$ from 
$(A_i,\psi _i)$ on one hand and product formulas for the Seiberg-Witten invariants on manifolds
that are  fiber sums of simpler manifolds. 
We begin by briefly (and with few details) recalling the scenario of the latter. 

Let $X_i$, $i=0,1$ be two compact smooth 4-manifolds and $\Sigma _i \hookrightarrow X_i$
embedded surfaces of the same genus and with 
$\Sigma _0 \cdot \Sigma_0 = - \Sigma _1 \cdot \Sigma _1$. In this setup one can construct the
fiber sum 
$$X=X_0 \# _{\Sigma _i} X_1$$ 
by cutting out tubular neighborhoods $N(\Sigma _i )$ in $X_i$ and gluing the manifolds
$X_i'=\overline{X_i \backslash  N(\Sigma _i )}$ along their diffeomorphic boundaries. 

Under certain conditions  one can calculate some of the Seiberg-Witten invariants of $X$ in
terms of the Seiberg-Witten invariants of the building blocks $X_i$ 
(see e.g. \cite{kn:szabo}). One 
accomplishes this by showing that from solutions $(B_i, \Phi _i)$, $i=0,1$ on $X_i$ one can
construct a solution $(B,\Phi )$ on $X$  (this isn't possible for any pair of 
solutions $(B_i, \Phi _i)$ but the details are not relevant to the present discussion). This is done
by inserting a ``neck" of length $r\ge 1$ between the $X_i '$ so 
as to identify $X$ with 
$$ X = X_0' \cup \left(  [0,r]\times Y \right) \cup X_1'$$
with $Y=\partial N(\Sigma _0 ) \cong \partial N(\Sigma _1) $. A partition of unity $\{ \varphi _0 
,  \varphi _1 \}$ is chosen for each value of $r\ge 1$ subject 
to the conditions: 
\begin{align} 
\varphi _i = 1 &  \mbox{ on } X_i' \nonumber \\
\varphi _i = 0 &  \mbox{ outside of  } X_i' \cup [0,r] \times Y \nonumber \\
|\varphi _i'| \le \frac{C}{r} & \mbox{ on } [0,r] \times Y \nonumber
\end{align}
An approximation $\Phi '$ of $\Phi$ is then defined to be $\Phi ' = \varphi _0 \, \Phi _0  + \varphi
_1 \, \Phi _1 $ (similarly for $B'$, a first
approximation for $B$). The measure of the  failure of  $(B' ,\Phi ')$  to solve the Seiberg-Witten
equations can be made as small as desired by making $r$ large.
The honest solution $(B, \Phi )$ is then sought in the form $(B',\Phi ') + (b,\phi) $ with $(b,\phi)
$ small. The correction term $(b,\phi )$ is found as a fixed point of the 
map 
$$(b,\phi ) \mapsto Z(b,\phi ) = -L_{(B', \Phi ')}^{-1}  \left( Q(b,\phi ) + \mbox{err} \right) $$
Here ``err" is the size of $SW(B', \Phi ')$ and $L$ and $Q$ are as in the previous section.
Choosing
$r$ large enough and $||(b,\phi )||$ small enough makes $Z$ 
a contraction mapping and so the familiar fixed point theorem for Banach spaces guarantees the
existence of a unique fixed point.

In the case of fiber sums there are product formulas that allow one to calculate the
Seiberg-Witten invariants of $X$ in terms of the 
invariants of the manifolds $X_i$. The formulas typically have the  form: 
\begin{equation}
SW_X (W_E) = \sum _{E_0 + E_1=E} SW_{X_0} (W_{E_0}) \cdot SW_{X_1} (W_{E_1}) 
\label{eq:prodformula}
\end{equation}
Due to the similarity of our construction of grafting monopoles to the one used to construct $(B,
\Phi )$ from $(B_i, \Phi _i)$, 
it is natural to ask if such or similar formulas exist for the present case, that is, can one calculate
$SW_X(W_{E_0 \otimes E_1})$ in terms of $SW_X(W_{E_0})$ and  
$SW_X(W_{E_1})$? The author doesn't know the answer. However, if they do exist, they can't
be expected to be as simple as \eqref{eq:prodformula}. 
 The reason for this can be understood by trying to take the analogy between our setup and that
for fiber sums further. 

In the case of fiber sums, once one has established that the two solutions $(B_i, \Phi _i )$ on
$X_i$ can be used to construct a solution $(B,\Phi )$ on $X$, 
one  needs to establish a converse of sorts. That is, one needs to show that every solution
$(B,\Phi )$ on $X$ is of that form. It is at this point where 
 the analogy between 
the two situations breaks down. It is conceivable  in our setup,
that there will be solutions for the \spin $(E_0\otimes E_1)\otimes W^+_0$ that 
can not be obtained as products of solutions for the \spins $E_i\otimes W^+_0$. Worse even,
there might be monopoles that 
can not be obtained as products of solutions for any \spins 
$F_j\otimes W^+_0$ with the  choice of  $F_j$, $j=0,1$ such that  $E=F_0\otimes F_1$ and $F_j\ne
0$. Those are the  monopoles where $\alpha ^{-1} (0)$ 
is connected. 
Thus if a product formula for our situation exists, it must in addition to a term
similar to  the 
right hand side of \eqref{eq:prodformula} also contain terms which count these
``undecomposable" solutions. But then again, they might not exist. 

The next section describes which solutions of the Seiberg-Witten equations for the \spin
$(E_0\otimes E_1)\otimes W^+_0$ are obtained as products of 
solutions for the \spins $E_i\otimes W^+_0$, $E=E_0\otimes E_1$. 
\section{The image of the multiplication map} \label{sec:5}
This section describes a partial converse to theorem \ref{theorem:main}. Recall that 
$$\Theta: {\cal M}^{SW}_X (W_E) \rightarrow {\cal M}^{Gr}_X (E) $$ 
is the map assigning a $J$-holomorphic curve to a Seiberg-Witten monopole. 
%
%
%
\begin{theorem} \label{theorem:converse}
Let $E=E_0\otimes E_1$ and let  $(A,\psi )$ be a solution of the 
Seiberg-Witten equations in the \spin $W_E$ with perturbation term $\mu = F_{A_0}^+ -
ir\omega /8$ and with $\psi = \sqrr \,(\alpha \otimes u_0 , \beta )$. 
Assume further that $J$ has been chosen from the Baire set ${\cal J}_0$ and that $\Theta (A,\psi )$ contains
no multiply covered components. If there exists an $r_0$ such  that for all $r\ge r_0$, 
$\alpha ^{-1}(0)$ splits into a disjoint union $\alpha ^{-1} (0) =\Sigma _0 \sqcup \Sigma _1$
with $[\Sigma _i]=$P.D.$(E_i)$ then $(A,\psi )$ lies in 
the image of the multiplication map 
$$ {\cal M}^{SW}_X(E_0) \times  {\cal M}^{SW}_X(E_1) \stackrel{\cdot}{\rightarrow}  {\cal
M}^{SW}_X(E_0\otimes E_1)$$
\end{theorem}

The proof of theorem \ref{theorem:converse} is divided into 3 sections. In section
\ref{section:one} we give the 
definition of $(A_i',\psi _i')$ - first approximations of Seiberg-Witten monopoles  $(A_i,\psi _i)$
for the \spin $W_{E_i}$ which when multiplied give the monopole 
$(A,\psi )$ from theorem \ref{theorem:converse}. Section \ref{section:two} shows that for large
values of $r$, $(A_i',\psi _i')$ come close to solving the
Seiberg-Witten equations. In the final section \ref{section:three} we show that $L_{(A_i',\psi
_i')}$ is surjective with inverse bounded independently of $r$. 
The contraction mapping principle is then used to deform the approximate 
solutions $(A_i',\psi _i')$ to honest solutions $(A_i,\psi _i)$. Section \ref{section:three} also
explains why $(A_0,\psi _0) \cdot (A_1,\psi _1) = (A,\psi )$.

We tacitly carry the assumptions of the theorem until the end of the section.
\subsection{Defining $(A_i',\psi _i')$} \label{section:one}

The basic idea behind the definition of $(A_i',\psi _i')$ is again that of grafting existing
solutions.
For example, one would like $\psi _0 '$ to be 
defined as the restriction of $\psi$ to a neighborhood of $\Sigma _0$ (under an appropriate
bundle isomorphism trivializing $E_1$ over that neighborhood) and 
to be the restriction of $\sqrr \, u_0$ outside that neighborhood. This is essentially how the construction
goes even though a bit more care is required, especially in 
splitting the connection $A$ into $A_0'$ and $A_1'$. 

To begin with, choose regular neighborhoods $V_0$ and  $V_1$ of $\Sigma _0$ and $\Sigma
_1$.  Once $r$ is large enough, these choices don't need 
to be readjusted for larger values of $r$. Choose, as in section \ref{section:four}, an open set
$U$ such that: 
\begin{align}
X&= V_0 \cup U \cup V_1 \nonumber \\
& U\cap \Sigma _i = \emptyset \nonumber 
\end{align}
Also, just as in section \ref{section:four}, 
arrange the choices so that $U\cap V_i$ contains a collar $\partial V_i \times [0,d]$ (with
$\partial V_i$ corresponding to $\partial V_i 
\times \{ d \}$) and choose $\delta > 0$ smaller than $d/1000$. Assume that the curves $\Sigma _i$ are connected, 
the general case goes through with little difficulty but with a bit more complexity of notation. 

Over $U\cup V_1$, choose a section $\gamma _0\in \Gamma (E_0;U\cup V_1)$ with $|\gamma
_0| = 1$.  Choose a connection $B_0$ on $E_0$ with respect to 
which $\gamma _0$ is covariantly constant over $U\cup V_1$, i.e. 
\begin{equation}
B_0 (\gamma _0 (x)) = 0 \quad \quad \forall x\in U\cup V_1
\end{equation}
Notice that such a connection is 
automatically flat over $U\cup V_1$. Choose a connection $B_1$ on $E_1$ such that
$B_0\otimes B_1 = A$ over $X$. Now define 
$\tilde{\alpha_1  }'\in \Gamma (E_1;U\cup V_1)$ and $\tilde{\beta _1}' \in  \Gamma
(E_1\otimes K^{-1};U\cup V_1)$ by:   
\begin{align}
\alpha = \gamma _0 \otimes \tilde{\alpha_1  }'  \\ 
\beta = \gamma_0 \otimes \tilde{\beta _1}' 
\end{align}
Proceed similarly over $V_0$. However, since some of the data is now already defined, more
caution is required. Choose a section 
$\gamma _1 \in \Gamma(E_1;V_0)$ with: 
\begin{equation} \begin{array}{ll}
 \gamma _1 =   \tilde{\alpha_1} ' & \mbox{ on } (U\cap V_0)\backslash (\partial V_0\times
[0,4\delta\rangle) \\
                                               & \\
  |\gamma _1 | = 1                     &  \mbox{ on } \left( V_0\backslash U\right)  \cup \left( \partial
V_0 \times [0,2\delta\rangle \right)   \\
\end{array}
\end{equation}
We continue by defining $\tilde{\alpha _0}'$ and $\tilde{\beta _0}'$ over $V_0$ by: 
\begin{align}
\alpha = \tilde{\alpha_0  }'  \otimes \gamma _1  \\ 
\beta =  \tilde{\beta _0}' \otimes  \gamma_1
\end{align}
Choose one forms $a_0$ and $a_1$ such that over $V_0$ the following two relations hold:
\begin{align}
(B_1 + i \, a_1) \, \gamma _1 & = 0  \\ 
(B_0 + i \, a_0)\otimes (B_1 + i \, a_1)  & = A 
\end{align}
With these preliminaries in place, we are now ready to define $(A_i ',  \psi _i ')$: 
\begin{equation}
\begin{array}{ll}
\tilde{\alpha_0} = \chi _{4\delta ,0} \, \tilde{\alpha_0  }' + (1- \chi _{4\delta ,0})  \gamma _0
\phantom{mmm} & \tilde{\beta_0} = \chi _{4\delta , 0} \, \tilde{\beta_0} '  \\ 
\tilde{\alpha_1} = (1-\chi _{4\delta ,0}) \, \tilde{\alpha_1  }' + \chi _{4\delta ,0}  \gamma _1 &
\tilde{\beta_1} = (1-\chi _{4\delta , 0}) \, \tilde{\beta_1} '  \\ 
A_0' = B_0 + i \chi _{4\delta ,0} \, a_0 & A_1' = B_1 + i \chi _{4\delta ,0} \, a_1\\ & \\
\end{array}
\end{equation}
\begin{lemma}  \label{lemma:props} The $(A_i',\psi _i ' )$ defined above, satisfy the following
properties:
\begin{itemize}
\item[\rm (a)] $A_0'\otimes A_1 ' = A$ on all of $X$.
\item[\rm (b)] $\tilde{\alpha _0} = \gamma _0$ on $(U\cap V_0) \backslash (\partial V_0 \times
[0,4\delta \rangle) $.
\item[\rm (c)] $F_{B_0} = 0 $ on $U\cup V_1$ and $F_{B_1+i\tilde{a_1}} = 0 $ on $V_0$. 
\item[\rm (d)] On $(U\cap V_0) \backslash (\partial V_0 \times [0,4\delta \rangle) $, $|\tilde{a_i} |$
and $|d\, \tilde{a_i}|$ converge exponentially fast to zero 
as $r\rightarrow \infty$. 
\end{itemize}
\end{lemma}
\begin{proof} (a)\qua This is trivially true everywhere except possibly on the support of $d\chi
_{4\delta,0}$ which is contained in   $U\cap V_0$. 
However, on $U\cap V_0$ we have $A=B_0\otimes B_1$ and $A=(B_0+ia_0)\otimes
(B_1+ia_1)$ and thus $a_0 + a_1 = 0$. 
In particular, $A_0'\otimes A_1'= B_0\otimes B_1 + i \chi_{4\delta ,0} \, (a_0 + a_1)=
B_0\otimes B_1 = A$.

(b)\qua  Notice that on  $(U\cap V_0) \backslash (\partial V_0 \times [0,4\delta \rangle) $, $\gamma_1
= \tilde{\alpha _1'}$. Thus, $\alpha = \gamma _0 
\otimes \gamma _1$ and $\alpha = \tilde{\alpha _0 '} \otimes \gamma _1$ imply that $\gamma
_0 =  \tilde{\alpha _0 '}$. The claim now follows from the 
definition of $\alpha _0$.

(c)\qua    Follows from the fact that both connection annihilate nowhere vanishing sections on the said
regions.

(d)\qua  On $(U\cap V_0) \backslash (\partial V_0 \times [0,4\delta \rangle) $ we have  $\alpha =
\gamma_0 \otimes \gamma _1$ and $\nabla ^A = 
\nabla ^{B_0+ia_0} \otimes \nabla ^{B_1+ia_1}$. Also, recall that $\nabla ^{B_0} \gamma _0
= 0$ and $\nabla ^{B_1+ia_1} \gamma _1=0$. Thus: 
\begin{align} \nabla ^a \alpha = (\nabla ^{B_0+ia_0} \otimes \nabla ^{B_1+ia_1}) (\gamma _0
\otimes \gamma_1) = i \, a_0 \gamma _0 \otimes \gamma _1  \nonumber 
\end{align}
This equation yields: 
\begin{equation} \, |a_0| =  \frac{ |\nabla ^A \alpha | }{ |\alpha | } 
\end{equation}
The claim follows now for $a_0$ by evoking the bounds \eqref{eq:mainbounds}. The same result holds for $a_1$
by the proof of part (a) where it is shown that 
$a_0+a_1 = 0$ on $U \cap V_0$. The statement for $da_i$ follows from part (c), the equation 
$F_A = F_{B_0+ia_0} + F_{B_1+i a_1}$ and the bounds \eqref{eq:mainbounds} for $|F_A|$.   
\end{proof}
\subsection{Pointwise bounds on $SW(A_i',\psi _i')$} \label{section:two}

\begin{proposition} Let $(A_i',\psi_i ')$ be defined as above, then there exists a constant $C$ and an
$r_0\ge 1$ such that for all $r\ge r_0$ the inequality 
$$ |SW(A_i',\psi_i ')|_x \le \frac{C}{\sqrr}$$
holds for all $x\in X$.
\end{proposition}
\begin{proof} We calculate the size of the contribution of each of the three Seiberg-Witten
equations separately. The only nontrivial part of the 
calculation is in the region of $X$ which contains the support of $d\chi _{4\delta ,0}$ i.e.\ in
$\partial V_0\times [4\delta , 8 \delta]$. 
We will tacitly use the results of lemma \ref{lemma:props} in the calculations below. 
\eject
{\bf a)\qua The Dirac equation } 
\vskip2mm
To begin with, we calculate the expression $D_A ((\tilde{\alpha_0} \otimes u_0 + \tilde{\beta
_0} )\otimes \gamma _1) $ in two 
different ways. On one hand we have:
\begin{align} D_A ((\tilde{\alpha_0} \otimes u_0 + \tilde{\beta _0} )\otimes \gamma _1) = D_A
(\alpha + \chi _{4\delta ,0} \beta) = (1- \chi _{4\delta ,0}) D_A \alpha 
+ d\chi _{4\delta ,0} . \beta \nonumber
\end{align}
On the other hand we get:
\begin{align}  D_A ((\tilde{\alpha_0} \otimes u_0 & + \tilde{\beta _0} )\otimes \gamma _1) = \\ 
\nonumber
                   & = \gamma_1 \otimes D_{A_0'} (\tilde{\alpha _0}\otimes u_0 + \tilde{\beta _0}) 
+ e^i.(\tilde{\alpha _0}\otimes u_0 + \tilde{\beta _0}) \otimes A_1'(\gamma_1) \\ \nonumber
                   & = \gamma_1 \otimes D_{A_0'} (\tilde{\alpha _0}\otimes u_0 + \tilde{\beta _0}) 
+ i(\chi_{4\delta , 0} -1) a_1 \gamma_1 \nonumber
\end{align}
Equating the results of the two calculations we obtain:
\begin{align} |\alpha|\cdot |  D_{A_0'} (\tilde{\alpha _0}\otimes & u_0 +  \tilde{\beta _0}) | =
|\gamma_1 \otimes  D_{A_0'} (\tilde{\alpha _0}\otimes u_0 + \tilde{\beta _0}) |  \le  \nonumber
\\
\nonumber 
                 & \le C\,(  |\alpha | \, |a_1 | + |D_A \alpha| + |\beta| ) \le \frac{C}{\sqrr} \nonumber
\end{align}
Since over $\partial V_0 \times [4\delta , 8\delta ]$, $|\alpha | \rightarrow 1$ exponentially fast as
$r\rightarrow \infty$ we obtain that: 
\begin{equation}
| D_{A_0'} (\tilde{\alpha _0}\otimes  u_0 +  \tilde{\beta _0}) |  \le \frac{C}{\sqrr} 
\end{equation}
\vskip5mm
{\bf b)\qua The $(1,1)$-component of the curvature equation } 
\vskip2mm
Again, we only calculate for $x\in \partial V_0 \times [4\delta , 8\delta ]$: 
\begin{align} F^{(1,1)} _{A_0'} - F^{(1,1)} _{A_0} - & \frac{ir}{8} \left( |\tilde{\alpha _0} |^2
-1 - |\tilde{\beta _0}|^2 \right) \omega  = 
\chi_{4\delta , 0} \, (d a_0)^{(1,1)} + \frac{ir}{8}  |\tilde{\beta _0}|^2 \omega   \nonumber   \\ 
              & = \chi_{4\delta , 0} \, (d a_0)^{(1,1)} + \frac{ir}{8 \,  |\alpha | ^2} |\chi _{4\delta , 0}
|^2\, |\beta |^2 \omega  \nonumber
\end{align}
Both terms in the last line converge in norm exponentially fast to zero on $\partial V_0 \times [4\delta , 8\delta ] $ as
$r\rightarrow \infty$. 
\vskip5mm
{\bf c)\qua The $(0,2)$-component of the curvature equation } 
\vskip2mm
Similar to the calculation for the $(1,1)$-component of the curvature equation on $\partial V_0
\times [4\delta , 8 \delta]$, we have for the 
$(0,2)$-component of the same equation: 
\begin{align} F^{(0,2)} _{A_0'} - F^{(0,2)} _{A_0} - & \frac{ir}{4}  \overline{\tilde{\alpha
_0}} \, \tilde{\beta _0} = \chi_{4\delta , 0} \, (d a_0)^{(0,2)} - 
\frac{ir}{4 \,  |\alpha | ^2} \chi _{4\delta , 0} \, \overline{\alpha} \, \beta  \nonumber
\end{align}
Once again, both terms on the right-hand side of the above equation converge in norm exponentially fast to zero as $r$
converges to infinity. The proofs for the case of $(A_1',\psi _1 ')$ 
are similar and are left to the reader. 
\end{proof}
\subsection{Surjectivity of $L_{(A_i', \psi _i ')}$ and deforming ${(A_i', \psi _i ')}$ to an exact
solution} \label{section:three}

The strategy employed here is very similar to the one used in section \ref{section:six} and we
only spell out part of the details. 
We start by showing that $L_{(A_0', \psi _0 ')}$ is surjective, the case 
$L_{(A_1', \psi _1 ')}$ is identical. 

We begin by asking ourselves when the equation 
\begin{equation}
L_{(A_0', \psi _0 ')} \xi _0 = g_0 \label{equation:main3} 
\end{equation}
has a solution $\xi _0 \in L^{1,2}(i\Lambda ^1 \oplus W^+_{E_0}) $ 
for a given $g_0\in L^2(i\Lambda ^0 \oplus i \Lambda ^{2,+} \oplus W^-_{E_0})$. 
Define the analogues of the isomorphisms $\Upsilon _i$ from section \ref{section:four} to be:
\begin{align} \Upsilon _0 : \mathbb{C}  \times & (U\cup V_1) \rightarrow \Gamma (E_0 ;
U\cup V_1) \mbox{ given by } \Upsilon _0 (\lambda , x) = 
\lambda \cdot \gamma _0 (x) \mbox{ and } \nonumber \\ 
 & \Upsilon _1 : \mathbb{C} \times V_0 \rightarrow \Gamma (E_1;V_0) \mbox { given by }
\Upsilon_1 (\lambda , x) = 
\lambda \cdot \gamma _1 (x)  \nonumber \nonumber  
\end{align}
Let $\gamma \in L^2 (i\Lambda ^0\oplus i\Lambda ^{2,+} \oplus W_0^-;U\cup V_1) $ be
determined by the equation 
$\chi_{25\delta , 0} \, g_0 = \Upsilon_0 (\gamma)$ on $U\cup V_1$
and $\varsigma\in L^2 ( i \Lambda ^0\oplus i\Lambda ^{2,+} \oplus W^-_E;V_0)$ be given by
the equation 
$\Upsilon _1^{-1} (\varsigma) = (1-\chi_{25\delta ,0}) \, g_0$ on $V_0$. Thus we can write
$g_0$ as: 
\begin{equation}
g_0 =  \chi _{100\delta , 0} \, \Upsilon_0 (\gamma) + (1- \chi _{4\delta , 0}) \, \Upsilon _1^{-1}
(\varsigma) \label{equation:uno}
\end{equation} 
This last form suggests that, in order to split equation \eqref{equation:main3} into two components
involving $L_{(A,\psi)}$ and $S$, one should search for $\xi _0$ in the form 
\begin{equation}
\xi _0 =  \chi _{100\delta , 0} \, \Upsilon_0 (\eta) + (1- \chi _{4\delta , 0}) \, \Upsilon _1^{-1}
(\kappa) \label{equation:dos}
\end{equation}
with $\eta \in  L^{1,2} (i\Lambda ^1 \otimes W_0^+;U\cup V_1) $ and $\kappa \in  L^{1,2} ( i
\Lambda ^1 \oplus W^+_E;V_0)$. 
Using relations \eqref{equation:uno} and \eqref{equation:dos} in \eqref{equation:main3}  one obtains
the analogue of equation \eqref{eq:split}:
\begin{align} \chi _{100\delta , 0} & \Upsilon _0  (T (\eta ) - \Upsilon _0 ^{-1} \Upsilon _1
^{-1} {\cal P} (d\, \chi _{4\delta , 0} , \kappa) - \gamma ) + \nonumber \\
 & +  (1-\chi _{4\delta , 0}) \Upsilon _1 ^{-1} (M(\kappa) + \Upsilon_1 \, \Upsilon _0 {\cal P}(d
\, \chi _{100 \delta ,0} , \eta) - \varsigma) = 0 \label{equation:long}
\end{align}
The operators $T'$ and $M'$ are defined over $U\cup V_1$ and $V_0$ respectively, through the
relations: 
\begin{align} L_{(A_0',\psi _0 ')} \, \Upsilon _0 &= \Upsilon _0  T \nonumber \\
 L_{(A_0',\psi _0 ')} \, \Upsilon _1^{-1}  &= \Upsilon _1^{-1}   M \nonumber
\end{align}
We use these operators, defined only over portions of $X$, to define the operators $T'$ and $M'$
defined on all of $X$ by: 
\begin{align} T' =& (1-\chi _{\delta ,0}) T + \chi_{\delta,0} S \nonumber \\
                   M' =& \chi_{200\delta,0} M + (1-\chi_{200\delta ,0}) L_{(A,\psi)} \nonumber 
\end{align}
Split equation \eqref{equation:long} into the following two equations:
\begin{align} T'(\eta) - \Upsilon _0 ^{-1} \Upsilon _1 ^{-1} {\cal P} (d\, \chi _{4\delta , 0} ,
\kappa) = \gamma \nonumber \\
M'(\kappa) + \Upsilon_1 \, \Upsilon _0 {\cal P}(d \, \chi _{100 \delta ,0} , \eta) = \varsigma 
\label{equation:longsplit}
\end{align} 
It is easy to see that solutions to the system of equations \eqref{equation:longsplit} provide
solutions to \eqref{equation:long} by multiplying the two lines 
with $\chi _{100\delta , 0} \Upsilon _0$ and $(1-\chi _{4\delta , 0} \Upsilon _1 ^{-1})$
respectively and adding them. 

The following lemma is the analogue of lemma 3.6, its proof is identical to that of lemma 3.6
and will be skipped here. 
\begin{lemma} For every $\epsilon > 0$ there exists an $r_\epsilon \ge 1$ such that for $r\ge
r_\epsilon$ the following hold:
\begin{align} ||(M' - L_{(A,\psi)})\,x || _2 \le & \epsilon\,  ||x||_2 \nonumber \\
       ||(T'-S)\,y||_2 \le & \epsilon \, ||y||_2 \nonumber 
\end{align}
Here $x\in L^{1,2}(i\Lambda ^1 \oplus W^+_E)$ and $y\in L^{1,2} (i\Lambda ^1 \oplus W^+_0
)$. 
\end{lemma}
The lemma allows us to replace the system \eqref{equation:longsplit}  by the system:
\begin{align} S(\eta) - \Upsilon _0 ^{-1} \Upsilon _1 ^{-1} {\cal P} (d\, \chi _{4\delta , 0} ,
\kappa) = \gamma \nonumber \\
L_{(A,\psi)}(\kappa) + \Upsilon_1 \, \Upsilon _0 {\cal P}(d \, \chi _{100 \delta ,0} , \eta) =
\varsigma  \label{equation:longsplit2}
\end{align} 
The process of solving \eqref{equation:longsplit2} is now step by step the analogue of solving \eqref{eq:splittttt}. 
In particular, we solve the first of the two equations in 
\eqref{equation:longsplit2} for $\eta$ in terms of $\kappa$: 
$$\eta = \eta (\kappa) = S^{-1} ( \Upsilon _0 ^{-1} \Upsilon _1 ^{-1} {\cal P} (d\, \chi _{4\delta
, 0} , \kappa ) + \gamma ) $$
Use this in the second equation of \eqref{equation:longsplit2} and rewrite it as: 
$$\kappa = L_{(A,\psi)}^{-1} \left(    \varsigma -     \Upsilon_1 \, \Upsilon _0 {\cal P}(d \, \chi
_{100 \delta ,0} , \eta(\kappa) )     \right) $$
To solve this last equation is the same as to find a fixed point of the map $Y:L^2(i\Lambda ^1
\oplus W^+_E) \rightarrow L^2(i\Lambda ^1 \oplus W^+_E) $ 
(the analogue of the map described by \eqref{eq:llast}) given by:
$$Y(\kappa)  = L_{(A,\psi)}^{-1} \left(    \varsigma -     \Upsilon_1 \, \Upsilon _0 {\cal P}(d \,
\chi _{100 \delta ,0} , \eta(\kappa) )     \right) $$
The proof of the existence of a unique fixed point of $Y$ follows from a word by word
analogue of the proof of theorem \ref{theorem:mainn} together with the discussion 
preceding the theorem. 

With the surjectivity of $L_{(A_i',\psi _i')}$ proved, the process of deforming $(A_i',\psi _i ')$
to an honest solution $(A_i,\psi _i)$ is accomplished
by the same method as used in section 3.4 and will be skipped here. 

To finish the proof theorem 5.1, we need to show that: 
$$(A_0, \psi _0) \cdot (A_1 , \psi _1) = (A,\psi)$$
This follows from the fact that as $r\rightarrow \infty$, 
the distance dist$\left( (A_i,\psi _i), (A_i'\psi _i ') \right)$ converges to zero,  together with the following
relations which follow directly from the 
definitions:
\begin{align} \tilde{\alpha _0}\otimes \tilde{\alpha _1} = & \,\alpha \nonumber \\
\tilde{\alpha _0} \otimes \tilde{\beta _1} + \tilde{\alpha _1} \otimes \tilde{\beta _0} = & \,\beta
\nonumber \\
A_0'\otimes A_1' = & \, A \nonumber 
\end{align}  

\Addresses\recd

\end{document}